\documentclass[12pt]{article}
\usepackage[a4paper, total={6in, 8in}]{geometry}
\usepackage{mypreamble}
\usepackage{authblk}

\DeclareSymbolFont{CMletters}{OML}{cmm}{m}{it}
\DeclareMathSymbol{\nu}{\mathord}{CMletters}{23}

\graphicspath{{./pngfigs/}}

\providecommand{\keywords}[1]{\textbf{\textit{Keywords---}} #1}

\title{Patch-Smoother and Multigrid \\for the Dual Formulation for Linear Elasticity} 
\author{Gabriele Rovi \thanks{gabriele.rovi@usi.ch} and Rolf Krause \thanks{rolf.krause@usi.ch}\\
Euler Institute, Universit{\aaa}  della Svizzera italiana, Lugano, Switzerland}

\begin{document}

\maketitle

\begin{abstract}
The dual formulation for linear elasticity, in contrast to the primal formulation, is not affected by locking, as it is based on
the stresses as main unknowns. Thus it is quite attractive for nearly incompressible and incompressible materials. 
Discretization with mixed finite elements will lead to 
--possibly large-- linear saddle point systems with a particular structure.

Whereas efficient multigrid methods exist for solving problems in mixed plane elasticity, to the knowledge of the authors, no multigrid methods are readily available for the general dual formulation. 

Two are the main challenges in constructing a multigrid method for the dual formulation for linear elasticity. First, in the incompressible limit, the matrix block related to the stress is semi-positive definite. Second, the stress belongs to $\bH_{\text{div}}$ and standard smoothers, working for $\bH^1$ regular problems, cannot be applied. 

We present a novel patch-based smoother for the dual formulation for linear elasticity. We discuss different types of local boundary conditions for the patch subproblems. Based on our patch-smoother, we build a multigrid method for the solution of the resulting saddle point problem and investigate its efficiency and robustness. Numerical experiments show that Robin conditions best fit the multigrid framework, leading eventually to multigrid performance.
\end{abstract}

\begin{sloppypar}

\section{Introduction}   
Computational Mechanics is of paramount importance in numerous engineering applications. From the physical point of view, the two main unknowns of interest are the displacement $\bu$ and the stress $\bsigma$. However, depending on the weak formulation, only one of the two is the unknown to be found, while the other one has to be post-processed.

The most known weak form is the primal formulation, which solves only for $\bu$, while $\bsigma$ must be post-processed. In this case, since $\bu \in \bH^1$, such formulation is straightforward to handle and can be easily generalized to non-linear contact problems. See for example \cite{Kik88}. For these cases, non-linear optimal solvers have already been proposed, see, e.g., \cite{KK01,BK12,KM11}. However, in many applications, the stress is a physical quantity of primary interest. Its post-processing by means of differentiation of the displacement would
reduce the order of approximation. Furthermore the energy functional of the primal formulation becomes unbounded for nearly incompressible or incompressible materials.

In contrast to the primal formulation, the dual formulation directly solves for $\bsigma$ and its functional remains bounded in the incompressible limit. The equilibrium condition and the symmetry of the stress are enforced by means of the displacement $\bu$ and the rotation $\btheta$, that take the role of Lagrange multipliers, see \cite{BF12}.  For the sake of simplicity, only homogeneous and isotropic materials will be examined. The mixed finite elements discretization gives rise to a --possibly large-- saddle point system. Here, we use the first order Raviart-Thomas elements for $\bsigma$, the discontinuous linear Lagrangian elements for $\bu$ and the skew-symmetric continuous linear Lagrangian tensors for $\btheta$. For further details on this choice, see again \cite{BF12}. 

Efficient multigrid methods exist for the saddle point system arising from the discretization of the dual formulation for plane linear elasticity. However, the finite element discretizations differ from the one proposed in \cite{BF12} which holds, in 2D and 3D, for nearly incompressible and incompressible materials. For example, in \cite{KS04} and in \cite{PW06}, for the discretization of the stress, the PEERS elements and the Arnold-Winther elements are respectively used, and only nearly incompressible materials are examined. To the authors' knowledge, an optimal solver for the discretization of the general dual formulation for linear elasticity for nearly-incompressible or incompressible materials has not been proposed yet. 

For the design of our multigrid method, two main challenges have to be taken into account. First, in the incompressible limit, the matrix block related to the stress is only semi-positive definite. Secondly, $\bsigma \in \textbf{H}_{\text{div}}$ and standard smoothers, working for $ \textbf{H}_{1}$ regular problems, cannot be applied. Therefore the patch smoother introduced by Arnold-Falk-Winther defined in \cite{AW97} for $\textbf{H}_{\text{div}}$ regular problems is extended to all the degrees of freedom related to $\bu$ and $\btheta$. Then, depending on the type of boundary conditions for the local problem, there are different possible scenarios: full Neumann, Dirichlet, or Robin conditions. We study the smoother for all these options. Numerical experiments show that the Robin conditions best fit the multigrid framework. Indeed, even in the case of aggressive coarsening, the complexity of the multigrid method can be optimal, if the Robin parameter is chosen carefully.

\section{Dual Formulation for Linear Elasticity}     
Let $\Omega \subset \mathbb{R}^d$, with boundary $\partial \Omega$, represent the solid body of interest. We can consider two disjoint open sets of the boundary, the Neumann boundary $\Gamma_N$ and the Dirichlet boundary $\Gamma_D$, satisfying $\overline{\partial \Omega}= \overline{\Gamma}_N \cup \overline{\Gamma}_D$, with $\Gamma_N \cap \Gamma_D= \emptyset$. The body is subject to the external volumetric force $\bff:\Omega \to \mathbb{R}^d$ on $\Omega$ and to the boundary force $\bg_N:\Gamma_N \to \mathbb{R}^d$ on $\Gamma_N$. On $\Gamma_D$ a prescribed displacement $\bg_D$ is enforced. We seek for sufficiently smooth displacements $\bu: \Omega \to \mathbb{R}^d$ and internal stresses $\bsigma: \Omega \to \mathbb{R}^{d,d}$ which solve the linear elasticity problem described in the following. \\
Linear elasticity requires to enforce equilibrium directly on the reference configuration $\Omega$, instead of its deformed configuration, which is unknown. Then the conservation of linear and angular momenta of the body $\Omega$ are given by the following relations:
\begin{subequations}
\begin{align}
\text{div} \bsigma = -\bff \qquad \text{in } \Omega \:,
\label{equilibrium_constraint}\\
\textbf{as}\: \bsigma =\textbf{0} \qquad \text{in } \Omega \:,
\label{symmetry_constraint}
\end{align}
\label{equilibrium_symmetry_constraint}
\end{subequations}
where we define
\begin{align}
\bas \bsigma: =\dfrac{1}{2}(\bsigma-\bsigma^T) \:,
\label{symmetry_of_stress}
\end{align} 
the antisymmetric part of $\bsigma$. Furthermore, the strains have to be small enough so that the kinematic relation between strains and displacements $\beps$  reduces to a linear one:
\begin{align}
\beps (\bu):= \dfrac{1}{2} \left( \nabla \bu + (\nabla \bu)^T \right) \:,
\label{eps_definition}
\end{align}
and the constitutive law has to be a linear relation between the stress and the strain:
\begin{align}
\bsigma= \pazocal{C} \beps.
\end{align}
In particular, for a homogeneous and isotropic body, the constitutive Hook's law reads as follows:
\begin{align}
\bsigma=\pazocal{C} \beps := 2\mu \beps + \lambda (\tr(\beps)) \bI \:,
\label{sigma_of_eps}
\end{align}
where $\pazocal{C}$ is the \textit{stiffness tensor}, $\mu$ and $\lambda$ are the Lam\eeee ${}$ parameters,  $\bI$ is the identity matrix in $d-$dimension and ${\tr:  \mathbb{R}^{d,d}\to  \mathbb{R}}$, defined such that ${ \text{tr} (\cdot) :=\sum_{i=1}^d [\cdot]_{ii}}$, is the trace operator. Since the stiffness tensor has full rank, we can also define its inverse, the so-called \textit{compliance tensor} $\pazocal{A} $, such that:
\begin{align}
\beps=\pazocal{A} \bsigma:=\dfrac{1}{2 \mu}\left( \bsigma- \dfrac{\lambda}{d \lambda + 2 \mu} (\tr \bsigma) \bI \right) \:.
\label{eps_of_sigma}
\end{align}
We can observe that, in the incompressible limit $\lambda \to \infty$, only the constitutive law expressed in terms of the compliance tensor is bounded and defined for any entry. Given any bounded $\bsigma$ and any bounded $\beps$ such that $ \tr{\beps}=0$:
\begin{subequations}
\begin{alignat}{4}
& \lim_{\lambda \to \infty} \pazocal{C} \beps= 2 \mu \beps \:,
\label{unboundness_of_C}\\
& \lim_{\lambda \to \infty} \pazocal{A} \bsigma=  \dfrac{1}{2 \mu}\left(\bsigma -\dfrac{\tr \bsigma}{d}\bI \right)  \:.
\label{boundness_of_A}
\end{alignat}
\end{subequations}
Finally we must also close the problem with the boundary conditions:
\begin{subequations}
\begin{align}
\bsigma \cdot \bn &=\bg_{N} \qquad \text{on } \Gamma_N \:,\\
 \bu&=\bg_D \qquad \text{on } \Gamma_D \:,
 \label{bd_for_strong_form}
\end{align}
\end{subequations}
where $\bn$ represents the outward normal of the boundary $\partial \Omega$. The strong formulation of linear elasticity for a homogenous and isotropic material is:
\begin{subequations}
\begin{alignat}{4}
 \text{div} \bsigma &= -\bff &&\quad \text{in } \Omega\:,
 \label{equilibrium_constraint}\\
\textbf{as}\: \bsigma& =\textbf{0} &&\quad \text{in } \Omega \:,
\label{symmetry_constraint}
\\
\pazocal{A} \bsigma&=\beps &&\quad \text{in } \Omega \:,\\
\bsigma \cdot \bn& =\bg_{N} &&\quad \text{on } \Gamma_N\:,\\
 \bu&=\bg_D &&\quad\text{on }  \Gamma_D\:.%
\end{alignat}
\label{strong_formulation}%
\end{subequations}
The primal weak form solves only for the displacement $\bu$, exploits the constitutive law (\ref{sigma_of_eps}) and is therefore unbounded in the incompressible limit for ${\lambda \to \infty}$. On the other hand, the dual formulation solves for the stress $\bsigma$, exploits (\ref{eps_of_sigma}) and is bounded even for ${\lambda \to \infty}$. In particular, the displacement $\bu$ and the rotation $\btheta$ can be interpreted as Lagrange multipliers for enforcing the equilibrium condition (\ref{equilibrium_constraint}) and the symmetry of the stress tensor (\ref{symmetry_constraint}). See \cite{BF12}, \cite{SSS11} for further references. We define:
\begin{subequations}
\begin{align}
\hspace{-2cm}
 H_{\text{div}}(\Omega)\ &:=\{ \bsigma \in \bL^2(\Omega): \text{div} \bsigma \in L^2(\Omega) \} \:,
 \label{hdiv_space} \\
\bSigma &:=\{ \bsigma \in \left[ H_{\text{div}}(\Omega)\right]^d\} \:,\label{Sigma_space} \\
\bSigma_{\bg_N} &:=\{ \bsigma \in \left[ H_{\text{div}}(\Omega)\right]^d: \bsigma \cdot \bn |_{\Gamma_N}=\bg_N \} \:,\label{Sigmagn_space} \\
\bSigma_{0} &:=\{ \bsigma \in \left[ H_{\text{div}}(\Omega)\right]^d: \bsigma \cdot \bn |_{\Gamma_N}=\textbf{0} \} \:, 
\label{Sigma0_space}\\
\bU&:=\left[L^2(\Omega)\right]^d \:, 
\label{V_space}\\
\bTheta&:=\{\bgamma \in \left[ L^2(\Omega)\right]^{d, d} : \bgamma+\bgamma^T=\textbf{0} \} \: .
\label{Theta_space}
\end{align} 
\end{subequations}
The dual weak formulation for the problem (\ref{strong_formulation}) can be written by using proper bilinear and linear forms. We seek for $\bsigma \in \bSigma_{\bg_N} $, $\bu \in \bU$ and $\btheta \in \bTheta$ such that:
\begin{subequations}
\begin{alignat}{3}
&a(\bsigma, \btau)+b(\btau,\bu)+c(\btau,\btheta)&&=f_a(\btau) && \qquad \forall \btau \in \bSigma_0 \:,\\
&b(\bsigma,\bv)&&=f_b(\bv)  && \qquad \forall \bv \in \bU\:,\\
&c(\bsigma,\bgamma)&&=0 &&\qquad  \forall   \bgamma \in \bTheta\:,
  \end{alignat}
  \label{weak_saddle_point}
  \end{subequations}
  where the bilinear and linear forms are defined as follows:
  \begin{subequations}
  \begin{align}
  a(\bsigma, \btau)&:=\int_{\Omega} \pazocal{A} \bsigma \cdot \btau\:,\\
  b(\btau,\bu)&:= \int_{\Omega}\bu  \cdot \text{div} \btau\:,\\
  c(\btau,\btheta)&:=\int_{\Omega} \btheta \cdot \bas \btau \:,\\
  f_a(\btau)&:=\int_{ \Gamma} (\btau\cdot \bn) \cdot \bg_{D} \:,\\
   f_b(\bv)&:=-\int_{ \Omega} \bff\cdot \bv  \:.
  \end{align}
  \label{weak_formual_formulation}%
  \end{subequations}
Since the problem is a saddle point, an LBB condition has to be fulfilled. In  \cite{BF12} it is proven that there exists $\beta>0$ such that:
 \begin{align}
 \inf_{
 \substack{\bu \in \bU\\
\btheta \in \bTheta
}
}
 \sup_{
 \substack{\bsigma \in \bSigma
}
}
\dfrac{b(\bsigma,\bu)+c(\bsigma,\btheta)}{\norm{\bsigma}_{\bSigma} (\norm{\bu}_{\bU}+\norm{\btheta}_{\bTheta})} \geq \beta >0
\:.
\label{dual_LBB} %
 \end{align}

\section{Mixed Discretization}
Let $\pazocal{T}_h=\{K_1,...,K_{N_e} \}$ be a shape-regular simplicial mesh of $\Omega$ with $N_e $ elements. The general simplex $K$ is a triangle, in 2D, or a tetrahedron, in 3D. The subscript $h$ represents the maximal diameter of $\pazocal{T}_h$. The continuous spaces $\bSigma$, $\bU$ and $\bTheta$ are discretized on the mesh $\pazocal{T}_h$ by means of the finite element method. We get $\bSigma_h$, $\bU_h$ and $\bTheta_h$. All the other known functions, such as $\bn$, $\bff$, $\bg_D$, $\bg_N$ need to be discretized as well. We obtain $\bn_h$, $\bff_h$, $\bg_{h,D}$, $\bg_{h,N}$. The linear forms in (\ref{weak_saddle_point}) are consequently approximated and their discretizations are denoted with the subscript $h$. The continuous LBB version (\ref{dual_LBB}) must be satisfied also in the discrete setting. Then there must exist $\beta_h >0$ such that
 \begin{align}
 \inf_{
 \substack{\bu_h \in \bU_h \\
\btheta_h  \in \bTheta_h 
}
}
 \sup_{
 \substack{\bsigma_h  \in \bSigma_h 
}
}
\dfrac{b_d(\bsigma_h,\bu_h  )+c_d(\bsigma_h,\btheta_h )}{\norm{\bsigma_h }_{\bSigma_h } (\norm{\bu_h }_{\bU_h }+\norm{\btheta_h }_{\bTheta_h })} \geq \beta_h  >0
\:.
\label{dual_discrete_LBB}
 \end{align}
To this purpose, we choose for the following triplet of the discrete spaces:
\begin{subequations}
\begin{align}
 \bSigma_h &= \left[RT_1(\pazocal{T}_h)\right]^d \:,\\
 \bU_h&=\left[DP_1(\pazocal{T}_h)\right]^d\:,\\
 \bTheta_h&=\left[P_1(\pazocal{T}_h) \right]^{d \times d} \cap \bTheta\:,
 \label{dual_choice_discrete_spaces}
 \end{align}
\end{subequations}
where by $RT_1$, $P_1$ and $DP_1$ we denote respectively the spaces of first order Raviart-Thomas, continuous linear Lagrangian and discontinuous linear Lagrangian functions. See \cite{BF12}, \cite{RKL09}. We seek for ${\bsigma_h \in \bSigma_{h,\bg_{h,N}}}$, ${\bu_h \in \bU_h}$, ${\btheta \in \bTheta_h}$ such that:
\begin{subequations}
\begin{alignat}{3}
&a(\bsigma_h, \btau_h)+b(\btau_h,\bu_h)+c(\btau_h,\btheta_h)&&=f_{h,a}(\btau_h) && \qquad \forall \btau_h \in \bSigma_{h,0} \:,\\
&b(\bsigma_h,\bv_h)&&=f_{h,b}(\bv_h)  && \qquad \forall \bv_h \in \bU_h\:,\\
&c(\bsigma_h,\bgamma_h)&&=0 &&\qquad  \forall   \bgamma_h \in \bTheta_h\:.
  \end{alignat}
  \label{discrete_weak_form_dual_formulation}%
  \end{subequations}
Once the spaces, their bases and degrees of freedom (dofs) are chosen, it is also possibile to reformulate (\ref{discrete_weak_form_dual_formulation}) in a vector-matrix form. Let $n$ and $m$ be the dimensions of the spaces respectively of the unknown and of the Lagrange multipliers, i.e. $n=\text{dim}(\bSigma_h)$ and ${m=\text{dim}(\bU_h)+\text{dim}(\bTheta_h)}$. Furthermore let $\by_h \in \bY_h=\mathbb{R}^n $ and $\bz_h\in \bZ_h=\mathbb{R}^m$ be the vectors that collect the values of the dofs respectively of $\bsigma_h$ and $[\bu_h,\btheta_h]^T$. The problem (\ref{discrete_weak_form_dual_formulation}) can be equivalently written in vector-matrix form as follows:
\begin{align}
\begin{bmatrix}
\bA_h & \bB_h^T\\
{\bB_h} & \textbf{0}\\
\end{bmatrix}
\begin{bmatrix}
\by_h\\
\bz_h
\end{bmatrix}
=
\begin{bmatrix}
\bff_h\\
\bh_h
\end{bmatrix}
\:.
\label{discrete_saddle_point_system}
\end{align}  
where ${\bA_h \in \mathbb{R}^{n \times n}}$,  ${\bB_h \in \mathbb{R}^{m \times n}}$. The residuals corresponding to the whole system, to its first and to its second components can be defined as well:
\begin{subequations}
\begin{align}
\br_{\:\:\:}
&:=
\begin{bmatrix}
\bff_h\\
\bh_h
\end{bmatrix}
-
\begin{bmatrix}
\bA_h & \bB_h^T\\
{\bB_h} & \textbf{0}\\
\end{bmatrix}
\begin{bmatrix}
\by_h\\
\bz_h
\end{bmatrix}
\:,
\label{discrete_saddle_point_system_residuals_full}
\\
\br_a&:=\bff_h-\bA_h \by_h- \bB_h^T \bz_h 
\:,
\label{discrete_saddle_point_system_residuals_a}
\\
\br_b&:=\bh_h-\bB_h \by_h
\label{discrete_saddle_point_system_residuals_a}
\:.
\end{align}
\label{discrete_saddle_point_system_residuals}%
\end{subequations}
\begin{remark}
\label{incompressibility_remark}
In the incompressible limit, the operator $\pazocal{A}$ has a non trivial kernel:
\begin{align}
\text{Ker}\left(\lim_{\lambda \to \infty }\pazocal{A} \right)&=\{\bsigma : \: \bsigma = \alpha \bI, \quad \alpha \in \mathbb{R}  \} \:,
\label{kernel_A_lambda_inf}
\end{align}
and the same happens to the bilinear form $a$ and to the discretized matrix $\bA_h$, that can be only symmetric semi-positive definite. The fulfillment of the discrete LBB conditions (\ref{dual_discrete_LBB}) makes the whole system (\ref{discrete_saddle_point_system}) solvable, but there is no guarantee on the invertibility of the single block $\bA_h$. 
\end{remark}
\begin{remark}
\label{FOSLS_prof_remark}
We observe that the space $\bTheta$ consists of $L^2$ skew-symmetric matrix-value functions. Therefore, in the discrete setting, for $d=2$ only one scalar function is required, while for $d=3$ three scalar functions are needed. Indeed we can write: 
 \begin{align}
  \bTheta_h&=\Bigg\{\begin{bmatrix}
  0 & -p\\
  p & 0
  \end{bmatrix}\:, \quad p \in P_1(\pazocal{T}_h) \Bigg\} && d=2\\
  \bTheta_h&=\Bigg\{\begin{bmatrix}
  0 & -p &-q\\
  p & 0 & -r\\
    q & r & 0
  \end{bmatrix} \:, \quad p,q,r \in P_1(\pazocal{T}_h)\Bigg\} && d=3\:.
\end{align}  
The formulation can be rewritten in terms of $p$ for $d=2$ or in terms of $p,q,r$ for $d=3$. Then, instead of the space ${\bTheta_h=\left[P_1(\pazocal{T}_h) \right]^{d \times d} \cap \bTheta}$, we use $\Theta_h=P_1(\pazocal{T}_h)$ or $\bTheta_h=\left[{P}_1(\pazocal{T}_h)\right]^3$. 
\end{remark}

\section{Multigrid and Patch Smoother}
\label{multigrid_section}
\begin{figure}[htbp!]
\subcaptionbox{2D node-related divergence-free function.}[0.49\linewidth]{
\includegraphics[width=0.3\textwidth]{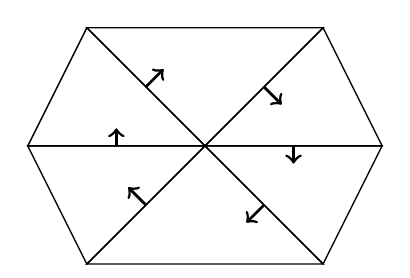}
}
\subcaptionbox{ 3D edge-related divergence-free function.
\label{divergencefree_functions_3D}
}[0.48\linewidth]{
\includegraphics[width=0.3\textwidth]{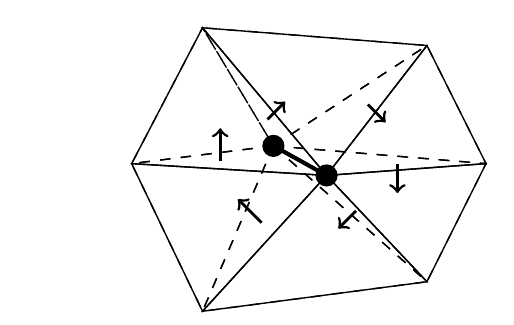}

}
\caption{Piecewise linear divergence-free functions.}
\label{divergencefree_functions}
\end{figure}

\begin{figure}[htbp!]
\subcaptionbox{All internal dofs.
\label{RT1_subspaces_a}
}[.4\linewidth]{
 \scalebox{.8}{
 \includegraphics[width=0.3\textwidth]{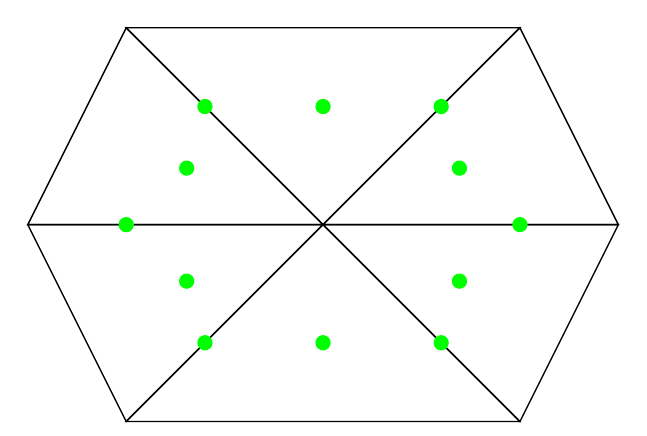}
 }
}\subcaptionbox{All dofs.
\label{RT1_subspaces_b}
}[.4\linewidth]{
 \scalebox{.8}{
\includegraphics[width=0.3\textwidth]{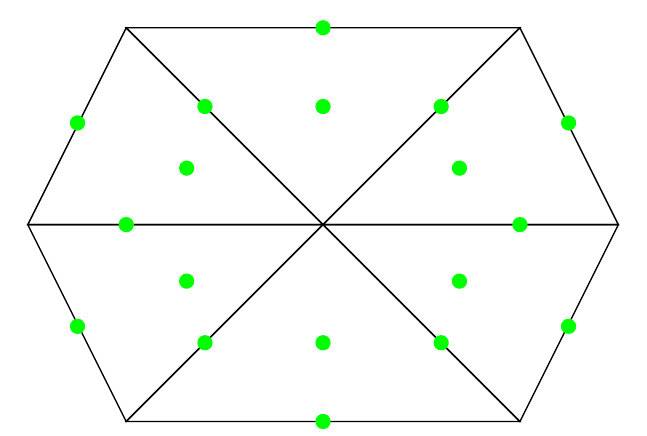}

 }
 }

 \caption{The green circles refer to $RT_1$ dofs. For simplicity, we count 2 dofs per each green circle, since we do have 2 dofs per face and 2 internal dofs.} 
\label{RT1_subspaces}
\end{figure}

\begin{figure}[htbp!]

\subcaptionbox{Only internal $\bsigma_h$ dofs.
\label{monolithic_subspace_a}
}[.4\linewidth]{
 \scalebox{.8}{
\includegraphics[width=0.3\textwidth]{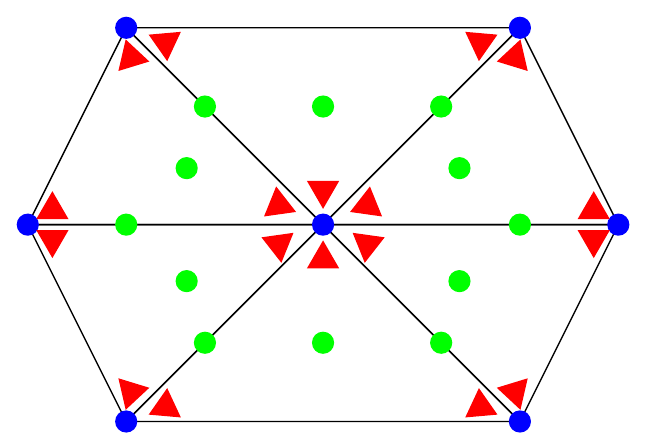}
}
}\subcaptionbox{All $\bsigma_h$ dofs.
\label{monolithic_subspace_b}
}[.4\linewidth]{
 \scalebox{.8}{
\includegraphics[width=0.3\textwidth]{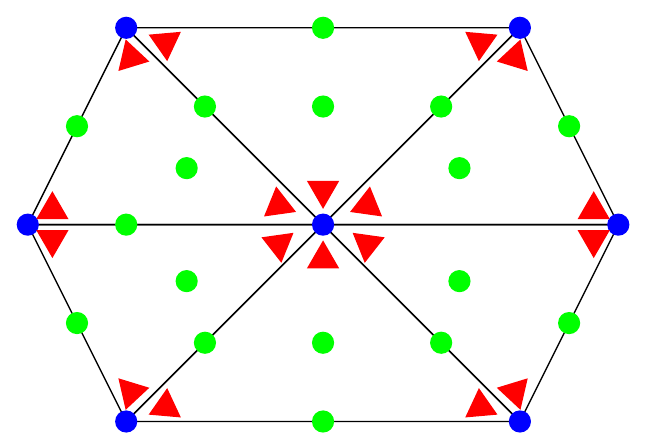}
}
} \includegraphics[width=0.15\textwidth]{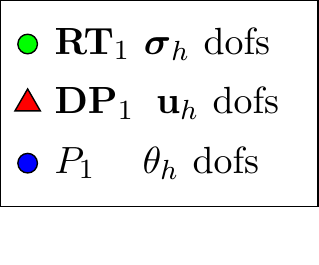}
\caption{Dofs of the subspace for the extended monolithic Arnold-Falk-Winther smoother. } 
\label{monolithic_subspace}
\end{figure}

A smoother is an iterative method that can rapidly damp the high-frequency components of the error. After the high frequencies are removed from the error, the smoother has a small impact on its low-frequency components. Furthermore, the larger the problem becomes, the more this behavior becomes evident and its convergence property deteriorates. The multigrid idea is to represent the error on coarser subspaces so that its low-frequency components become high-frequency components on these coarser subspaces and can still be easily damped with proper smoothers. In this way, optimal convergence is achieved: the number of iterations is independent of the dimension of the problem and of the number of levels used. See \cite{Xu92}, \cite{BY10}, \cite{Xu96}, \cite{CZ94}, \cite{GM90}, \cite{Man84}. To reach this goal, two are the main ingredients: a hierarchy of nested spaces and a smoother for each level. \\
Let ${\{\pazocal{T}_0 \}_{j=0}^J}$ be  a sequence of nested triangulations such that ${\pazocal{T}_0 \subset \pazocal{T}_1 \subset \cdots \subset \pazocal{T}_{J-1} \subset \pazocal{T}_J:=\pazocal{T}_h}$. On these meshes, we can define a sequence of nested subspaces ${\bY_0 \subset \bY_1 \subset ... \subset \bY_{J-1} \subset\bY_J:=\bY_h}$, where each $\bY_j$ is the coarse space on the level $j$ of the space $\bY_J$, related to the mesh $\pazocal{T}_j$. We also denote by ${\bPi_j^{j
+1}: \bY_j \to \bY_{j+1}}$ the interpolation operator between the levels $j$ and $j+1$. 
In a similar way, we can also consider a sequence of nested subspaces ${\bZ_0 \subset \bZ_1 \subset ... \subset \bZ_{J-1} \subset\bZ_J=\bZ_h}$. Let ${\bQ_j^{j
+1}: \bZ_j \to \bZ_{j+1}}$ be the interpolation operator between the levels $j$ and $j+1$. For $j=J-1,\hdots,0$ we get the coarse problems for the computation of the coarse correction $[\by_j,\bz_j]^T$:
\begin{align}
\begin{bmatrix}
\bA_j & \bB_j^T\\
{\bB_j} & \textbf{0}\\
\end{bmatrix}
\begin{bmatrix}
\by_j\\
\bz_j
\end{bmatrix}
=
\begin{bmatrix}
\bff_j\\
\bh_j
\end{bmatrix}
\:,
\label{discrete_coarse_problem}
\end{align}  
where
\begin{align}
\bA_j &:=\left[ \bPi_j^{j+1}\right]^T \bA_{j+1} \bPi_j^{j+1}\:,
\label{discrete_coarse_problem_A}
\\
\bB_j&: = \left[ \bQ_j^{j+1}\right]^T \bB_{j+1} \bPi_j^{j+1}\:,
\label{discrete_coarse_problem_B}\\
\bff_j &:= \left[ \bPi_j^{j+1}\right]^T \left( \bff_{j+1}- \bA_{j+1} \by_{j+1} \right)\:,
\label{discrete_coarse_problem_f}\\
\bh_j&:= \left[ \bQ_j^{j+1}\right]^T (\bh_{j+1}-\bB_{j+1} \by_{j+1}  )\:.
\label{discrete_coarse_problem_h}
\end{align}
Corrections on the fine and on all the other coarse levels can be computed by using a proper smoother. Standard smoothers like the component-wise Gauss-Seidel method and the conjugate gradient method cannot be used in this context, because the system is a saddle point. To overcome the saddle point structure, the Schur complement method or Uzawa's method, which requires to invert the single block $\bA_h$, could be considered. For a discussion on these methods, see for example \cite{Bra07}. Nevertheless, the Remark \ref{incompressibility_remark} implies that this direction cannot be examined and a monolithic approach is fundamental. 
Furthermore, the smoother to be used must be able to damp the divergence-free error components which show up since $ \bsigma \in \bH_{\text{div}}$. In order to tackle the divergence-free components of the error, two main strategies have been proposed in the existing literature. The first one is based on the Helmholtz decomposition and has been developed by \cite{Hip97}, \cite{HT00}, \cite{HX07}, \cite{XCN09}, \cite{KV12} and in other references therein. The second one has been examined by Arnold-Falk-Winther and aims to directly capture the divergence-free components of the error. See \cite{Arn},  \cite{Arn08}, \cite{AW97},\cite{AFW00}. The latter option will be investigated, since it can be easily extended to a local full monolithic approach, so that, in the incompressible limit, the relation (\ref{kernel_A_lambda_inf}) is not a problem anymore.\\
The Arnold-Falk-Winther's smoother is a sequential subspace correction method that tackles directly divergence-free components not visible from coarser meshes. To this aim, in contrast to the component-wise Gauss-Seidel, local corrections are not computed sequentially on one-dimensional subspaces, but on larger subspaces. In particular, the subspaces are defined on patches. A patch $\pazocal{P}_n$ related to a node $n$ is the set of all elements which share that node. This kind of patch can be used for any dimension $d$. In 3D we can also consider a patch $\pazocal{P}_{n,m}$ which is the set of all elements sharing the edge connecting the nodes $m$ and $n$. A representation of these patches and of piecewise linear divergence-free components is given in Figure \ref{divergencefree_functions}. However, from now on, we will focus on the 2D case. Concerning $\text{RT}_{1}$ functions, the subspace can be defined as the set of all internal dofs. See Figure \ref{RT1_subspaces_a}. But it can also be enlarged to all boundary dofs, as in Figure \ref{RT1_subspaces_b}.
${}$\\
The local subspaces for the extended monolithic Arnold-Falk-Winther's smoother must take into account also the dofs of the Lagrange multipliers. Therefore the subspaces of Figure \ref{RT1_subspaces} are extended as in Figure \ref{monolithic_subspace}. 
${}$\\
To this aim, for each level $j$, we define the decomposition ${\bY_j =\sum_{i=1}^{n_j} \bY_{j,i}}$, where each of the $n_j$ subspaces $ \bY_{j,i}$ is related to $\bY_j $ by means of the  interpolation operator $\bPi_{j,i} : \bY_{j,i} \to \bY_j $. Similarly we do for ${\bZ_j=\sum_{i=1}^{n_j} \bZ_{j,i}} $, with interpolation $\bQ_{j,i} : \bZ_{j,i} \to \bZ_j $. The local problem on the $i-$th subspace of level $j$, whose dofs are represented in Figure \ref{monolithic_subspace}, is the following:
\begin{align}
\begin{bmatrix}
\bA_{j,i} & \bB_{j,i}^T\\
{\bB_{j,i}} & \textbf{0}\\
\end{bmatrix}
\begin{bmatrix}
\by_{j,i}\\
\bz_{j,i}
\end{bmatrix}
=
\begin{bmatrix}
\bff_{j,i}\\
\bh_{j,i}
\end{bmatrix}
\:,
\label{discrete_ji_saddle_point_system}
\end{align}
where
\begin{align}
\bA_{j,i} &:=\bPi_{j,i}^T \bA_{j} \bPi_{j,i}\:,
\label{discrete_ji_coarse_problem_A}
\\
\bB_{j,i}&: =\bQ_{j,i}^T \bB_{j} \bPi_{j,i}\:,\\
\bff_{j,i} &:=\bPi_{j,i}^T \left( \bff_{j}- \bA_{j} \by_{j} \right)\:,\\
\bh_{j,i}&:=\bQ_{j,i}^T (\bh_{j}-\bB_{j} \by_{j}  )\:.
\label{discrete_ji_coarse_problem}
\end{align}
However, we are still free to choose how to deal with the local boundary conditions for the stress. In particular, full Neumann and full Dirichlet boundary conditions are represented for the local problem respectively in Figure \ref{monolithic_subspace_a} and in Figure \ref{monolithic_subspace_b}. In the following the full Neumann, Dirichlet, and Robin approaches will be examined.

\section{Full Neumann boundary conditions}
\label{arnold_smoother_dual_section}
\begin{figure}[htbp!]
\centering
\subcaptionbox{\small Starting patch.
\label{enlarged_patch_nu_a}
}[.25\linewidth]{
 \includegraphics[width=0.2\textwidth]{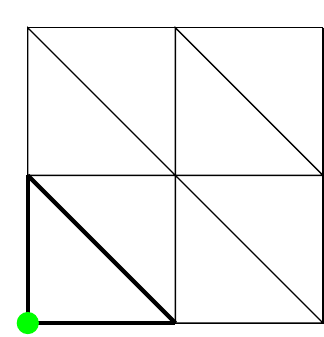}
  }
  \subcaptionbox{\small Enlarged patch.
\label{enlarged_patch_nu_b}
}[.25\linewidth]{
 \includegraphics[width=0.2\textwidth]{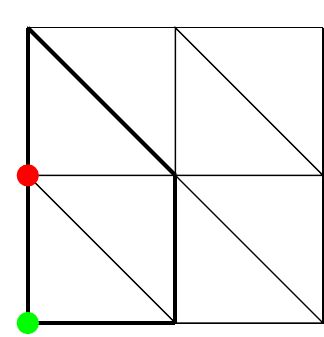}
  } \subcaptionbox{\small Enlarged patch.
\label{enlarged_patch_nu_c}
}[.25\linewidth]{
 \includegraphics[width=0.2\textwidth]{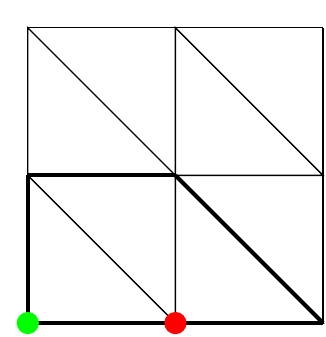}
 }
  \includegraphics[width=0.1\textwidth]{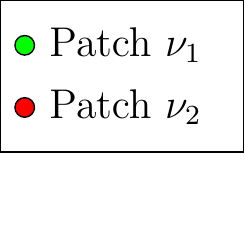}
 \caption{Starting from the patch built on the node $\nu_1$ of Figure \ref{enlarged_patch_nu_a}, we can have two possible enlarged patches, Figure \ref{enlarged_patch_nu_b} or Figure \ref{enlarged_patch_nu_c}, depending on the choice of $\nu_2$.}
\label{enlarged_patch}
\end{figure}

Neglecting the $\bsigma_h$ dofs on the boundary, as in Figure \ref{monolithic_subspace_a}, is equivalent to enforce local homogeneous Neumann boundary conditions. In linear elasticity, a body, that on the boundary is subject only to external forces, can freely translate and rotate. For $d=2$ the rigid body motions are three in total, given by two translations and one rotation, while for $d=3$ there are three rotations and three translations. Focusing on the case $d=2$ with only Neumann boundary conditions, we would expect a local system that can be solved only up to the rigid body motions. In the saddle point problem (\ref{discrete_ji_saddle_point_system}) related to the subspace of Figure \ref{monolithic_subspace_a}, such property shows up in the linear dependency of the constraints. Indeed the local matrix ${\bB_{j,i} \in \mathbb{R}^{m_{j,i} \times n_{j,i}}}$ of the problem (\ref{discrete_ji_saddle_point_system}) has not full rank if ${m_{j,i}>n_{j,i}}$. In particular, if free rigid body motions are permitted, it should happen that ${m_{j,i}=n_{j,i}+3}$. However, this is not always the case. 
In Figure \ref{table_1_2_3_elements_patch_subspace}, we notice that for a single element, $n_{j,i}=4$, $m_{j,i}=9$ and $m_{j,i}-n_{j,i}=5$. For two elements, $ n_{j,i}=12$, $m_{j,i}=16$ and $m_{j,i}-n_{j,i}=4$. Finally, for three elements $ n_{j,i}=20$, $m_{j,i}=23$ and $m_{j,i}-n_{j,i}=3$. Thus a patch of the kind of Figure \ref{monolithic_subspace_a} needs to be built at least on three elements. Otherwise the local matrix $\bB_{j,i}$ has too many linear dependent rows, even if we accept free rigid body motions.
\newcommand{\sr}{\rule[-0.45cm]{0pt}{3cm}}
\begin{figure}[htbp!]
\begin{tabular}{|c|c|c|c|c|c|}
\hline 
Legend dofs & Subspace &
\includegraphics[width=0.08\textwidth]{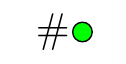}
&
\includegraphics[width=0.08\textwidth]{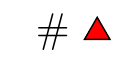}
&
\includegraphics[width=0.08\textwidth]{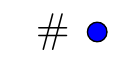}
&
\includegraphics[width=0.08\textwidth]{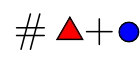}
\tabularnewline
\hline 
\multirow{3}{2cm}[0.cm]{
\vspace{-3cm}
\includegraphics[width=0.1\textwidth]{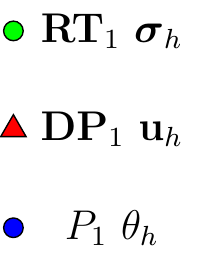}
} &\sr 

  \scalebox{.8}{
\includegraphics[width=0.35\textwidth]{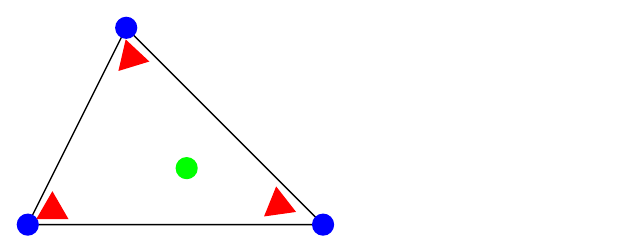}

}
& 
4 
&
6 
& 
3 
& 
9 

\tabularnewline
\cline{2-6} 
 & \sr 
   \scalebox{.8}{
\includegraphics[width=0.35\textwidth]{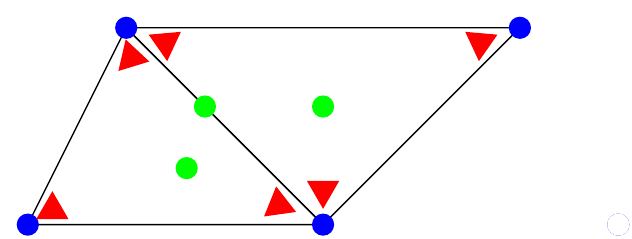}

}
& 
12 

&
12 

& 
4 

& 
16 

\tabularnewline
\cline{2-6} 
\cline{2-2} 
 & \sr
    \scalebox{.8}{
\includegraphics[width=0.35\textwidth]{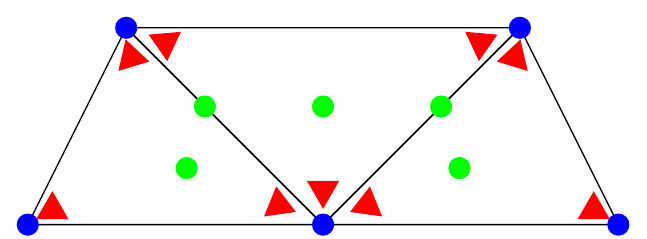}
 }
 
& 
20 

&
18 

& 
5 

& 
23 

\tabularnewline
\cline{2-6} 
\hline
 
\end{tabular}
\caption{The table shows how many dofs are related to the unknowns $\bsigma_h$ and how many to the Lagrange multipliers $\boldsymbol{u}_h$ and $\theta_h$. For full Neumann boundary conditions, only rigid body motions should be allowed, so the difference of Lagrange multiplier and stress dofs should be equal to 3. This happens if we consider a patch of at least 3 elements, where we do have $\# \text{dofs}(\boldsymbol{u}_h+\theta_h)=23$ and $\# \text{dofs}(\bsigma_h)=20$ and $\# \text{dofs}(\boldsymbol{u}_h+\theta_h)-\# \text{dofs}(\bsigma_h)=3$.}
\label{table_1_2_3_elements_patch_subspace}
\end{figure}
${}$\\
For fine enough meshes, the patches whose internal node is not on the boundary always have at least three triangles. However, on the border, no matter how fine the mesh is, a patch can consist of only one or two elements. Since this is the situation to be avoided, the definition of a  patch must be appropriately generalized. For a patch related to the node $\nu_1$, we do the following
\begin{itemize}
\item if it consists of at least three elements, then it is accepted;
\item otherwise, all the elements of the patch related to the node $\nu_2$ are added to the current patch, where $\nu_2 \neq \nu_1$ and $\nu_2$ belongs to one of the elements of the patch; if the new patch does not have at least three elements, the procedure is repeated for another node $\nu_3$ of the new patch, with ${\nu_3 \neq \nu_2, \nu_1}$; repeat until the patch consists of at least three elements; see Figure \ref{enlarged_patch} for an example of this procedure.
\end{itemize}

\def \ra{1.5}
\def \rb{3}

\def \rc{5}
\def \rd{10}

\begin{figure}[htbp!]
\centering
\hspace{-3cm}
 \includegraphics[width=0.6\textwidth]{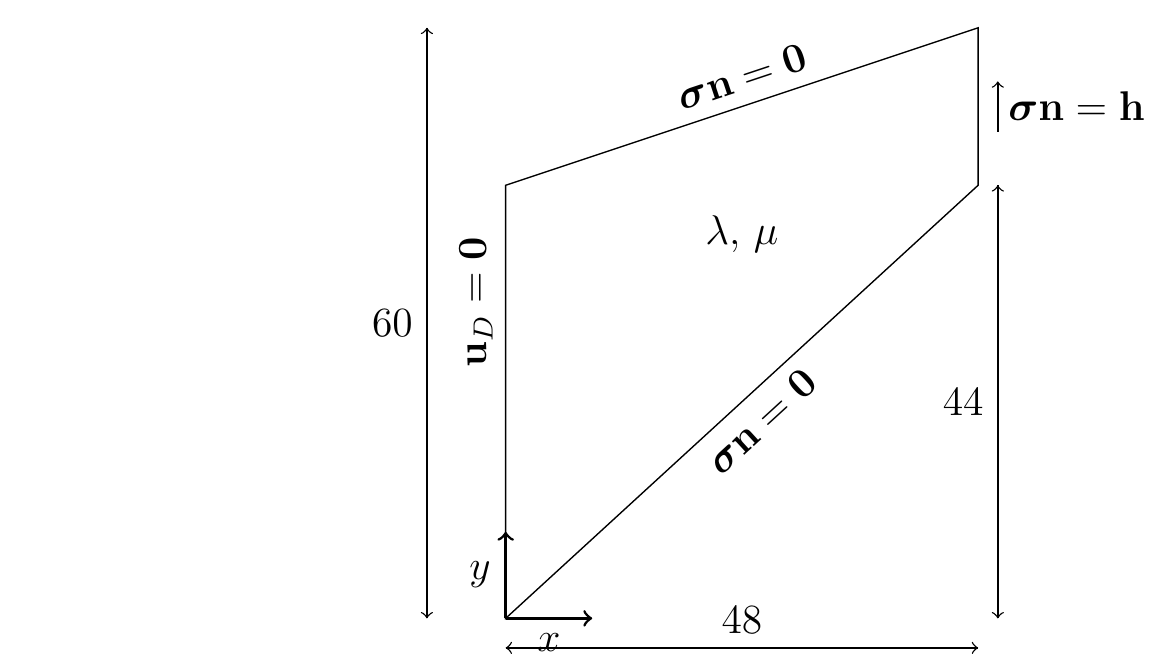}
\caption{Cook's membrane problem.}
\label{Cook_problem}
\end{figure}

\begin{figure}[htbp!]
    \scalebox{.7}{
{
 \includegraphics[width=0.6\textwidth]{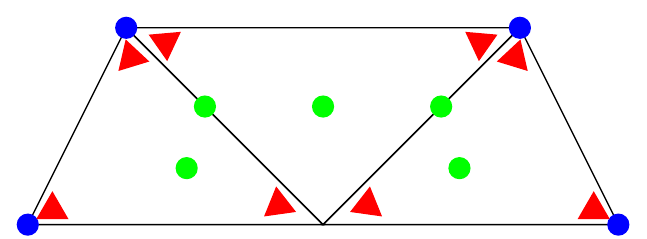}

}
 }  
 \scalebox{.7}{
{
 \includegraphics[width=0.6\textwidth]{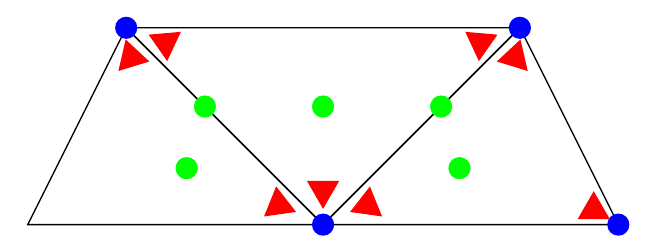}
}
 }
 \scalebox{.7}{
{
 \includegraphics[width=0.6\textwidth]{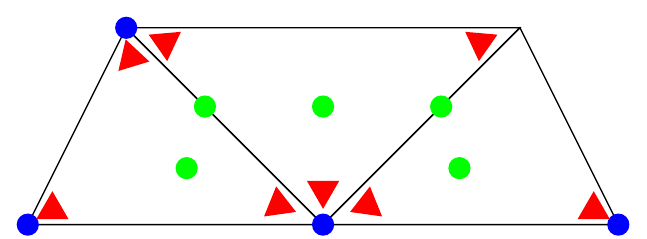}
}
 }
 \caption{Examples of patch built on the boundary where the displacement and rotation dofs related to a given node are removed from the subspace.}
 \label{patch_subspace_removal_1_rot_1_disp}
 \end{figure}

\begin{figure}[htbp!]
\centering
\includegraphics[width=0.263\textwidth]{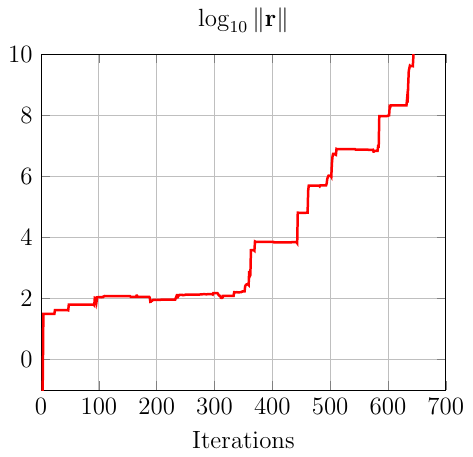}
\includegraphics[width=0.25\textwidth]{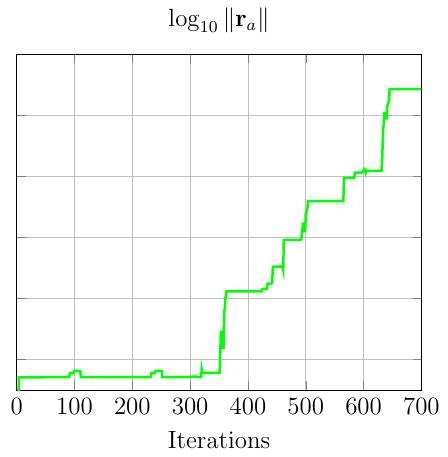}
\includegraphics[width=0.25\textwidth]{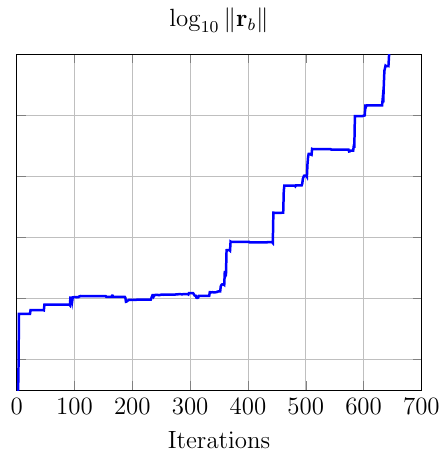}
\caption{$\log_{10}$ of the Euclidean norm of the residual for the Arnold-Falk-Winther smoother applied to dual formulation with subspaces of the type of Figure \ref{patch_subspace_removal_1_rot_1_disp} for the Cook's problem in Figure \ref{Cook_problem}.}
\label{patch_subspace_removal_1_rot_1_disp_result}
\end{figure}

The enlargement of the patches on the boundary can not be necessary in the case of the Dirichlet boundary. However, for simplicity, the process is carried out anyhow, independently of the type of boundary. Once the enlarged patches have been defined, the local problem can not still be solved, because rigid body motions are allowed. To overcome this difficulty, different are the strategies we could think of:
\begin{enumerate}
\item Removing rigid body motions from the local system, not considering the dofs related to one displacement, in both directions, and one rotation; of course, on a single patch, many combinations can be taken into consideration. 
\item Solving the system up to rigid body motions, searching for the local displacement and rotation which have zero-average on the patch.
\end{enumerate}
The first strategy is tested with the Cook's membrane problem, which is represented in Figure \ref{Cook_problem}. 
In figure \ref{patch_subspace_removal_1_rot_1_disp}, we represent the case in which, from the subspace of Figure \ref{monolithic_subspace_a}, a rotation and a displacement dofs are removed. Since on a given patch, many are the rotation and displacement dofs, as many combinations of removal of dofs can be made. For simplicity, after having locally numbered all the dofs, the first two components of the displacement and the last rotation of the patch are removed. The result of this specific case is represented in Figure \ref{patch_subspace_removal_1_rot_1_disp_result}. However, no matter how the removal process of the dofs is realized, a non-convergent smoother is recovered.

\begin{figure}[htbp!]
\centering
\includegraphics[width=0.265\textwidth]{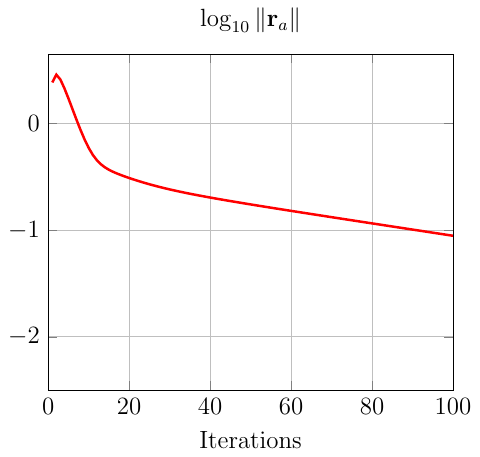}
\includegraphics[width=0.25\textwidth]{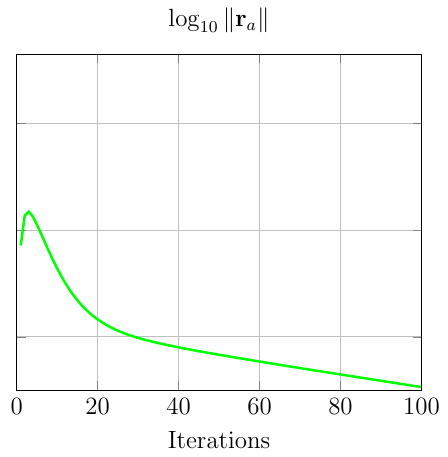}
\includegraphics[width=0.25\textwidth]{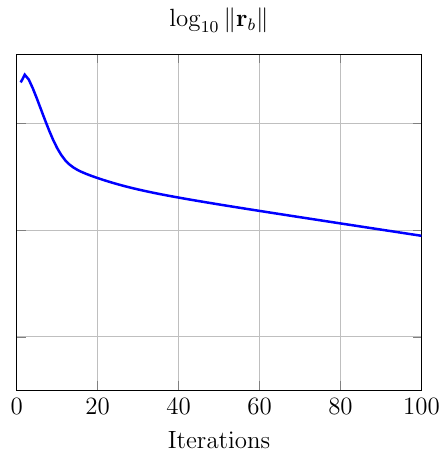}
\caption{$\log_{10}$ of the Euclidean norm of the residual for the Arnold-Falk-Winther smoother applied to the dual formulation for the Cook's problem in Figure \ref{Cook_problem}. Locally (\ref{zero_average_condition}) is enforced. Parameters: $\mu=1$, $\lambda=1$, $\text{number of smoothing steps}=100$. }
\label{zeroaveragemultipliersresiduals}
\end{figure}

The second strategy requires the local displacement and the local rotation to have zero average on the patch. For $d=2$:
\begin{align}
 \int _{\text{Patch}} u_{h,x} =0 \:,\qquad
 \int _{\text{Patch}} u_{h,y} =0 \:,\qquad
 \int _{\text{Patch}} \theta_h =0 \:,
 \label{zero_average_condition}
\end{align}
where we denote the first and second components of $\bu_h$ by $u_{h,x}$ and $u_{h,y}$, while the rotation $\theta_h$ is a scalar quantity. By applying the Arnold-Falk-Winther smoother which locally enforces (\ref{zero_average_condition}) to the Cook's problem of Figure \ref{Cook_problem}, we obtain the results in Figure \ref{zeroaveragemultipliersresiduals}. In the first iterations, the norm of the residual increases. Then it starts to decrease relatively fast and, after some iterations, very slowly. It is clear that the norm $\norm{\br}$ is governed by $\norm{\br_b}$,  since  $\norm{\br_b} \gg \norm{\br_a}$. So, by adding the condition (\ref{zero_average_condition}) the system is now solvable, but the constraints are not properly captured. It is not true that locally the solution must satisfy the zero-average condition for the displacement and the rotation. This is just a trick to make the local system solvable, but it also suggests that the choice of the subspace is not optimal. Besides, the most important iterations of a smoother are the initial ones. But here we see that many iterations are required to have a decrease in the residual. For this reason, another strategy is necessary.

\section{Full Dirichlet and Robin boundary conditions}
\label{patch_enriched_with_stress_dofs_on_the_boundary_section}
\begin{figure}[htbp!]
\centering
\includegraphics[width=0.6\textwidth]{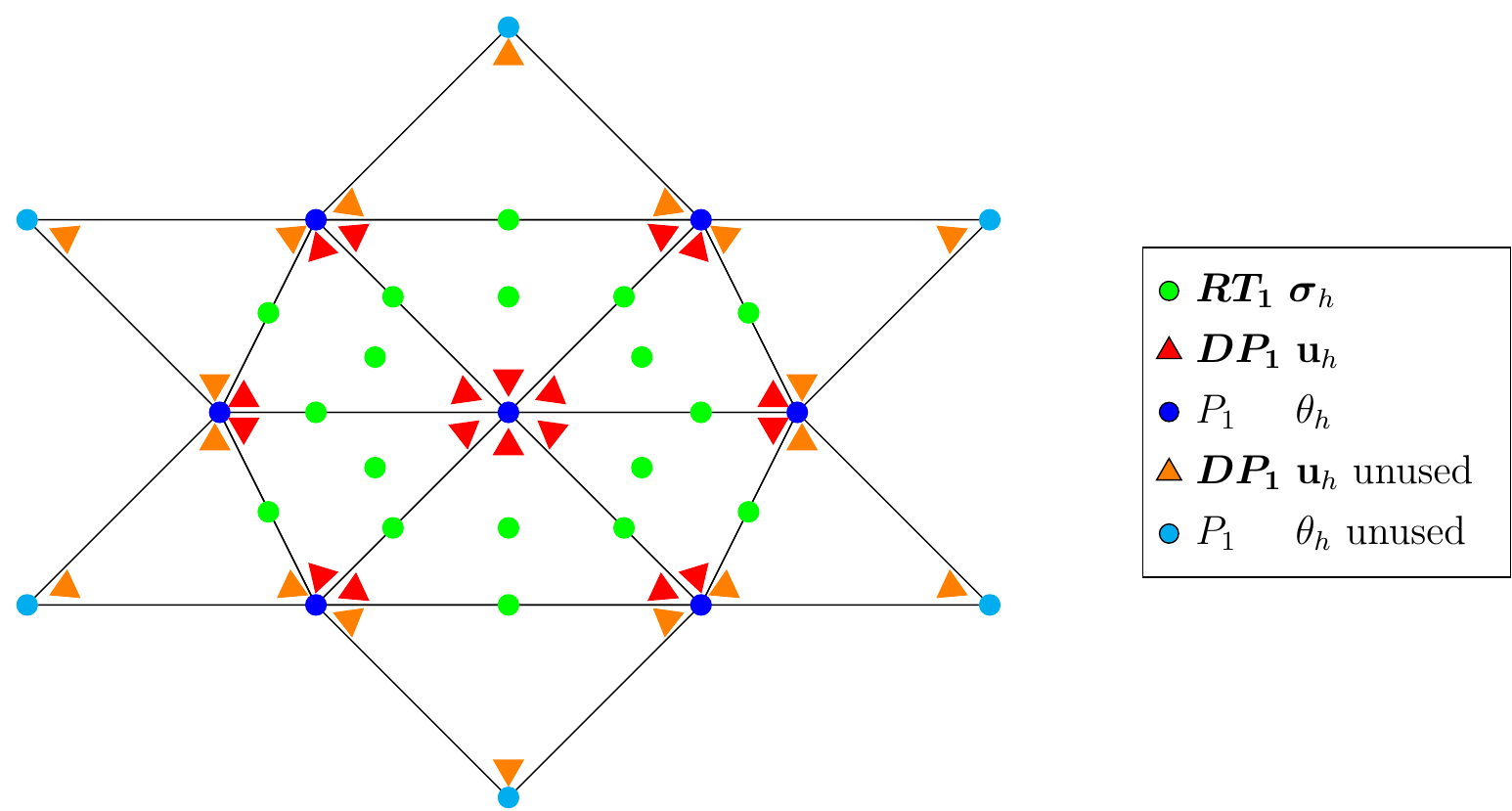}
\caption{Subspace referred to a patch considering $\boldsymbol{RT_1}$ dofs on the border. All the constraints outside the patch, related to the orange and cyan dofs, are unused but should somehow influence the local solution, since they communicate with the $\boldsymbol{RT_1}$ dofs on the border. }
 \label{RT1PatchExternalDofs}
\end{figure}

	\begin{figure}[htbp!]
	\centering
\includegraphics[width=0.275\textwidth]{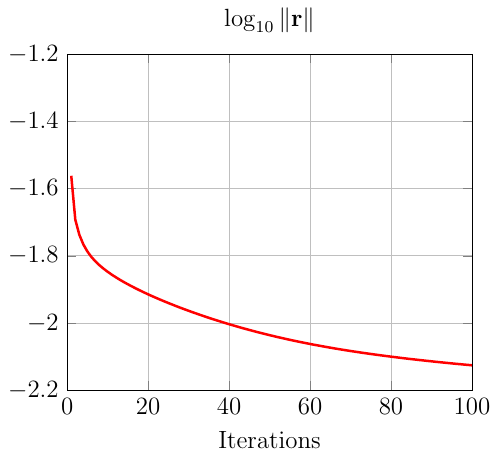}	
\includegraphics[width=0.25\textwidth]{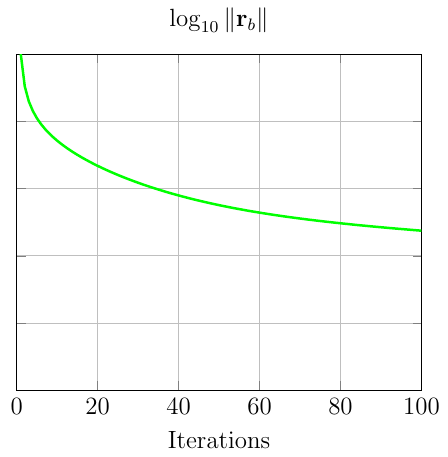}	
\includegraphics[width=0.25\textwidth]{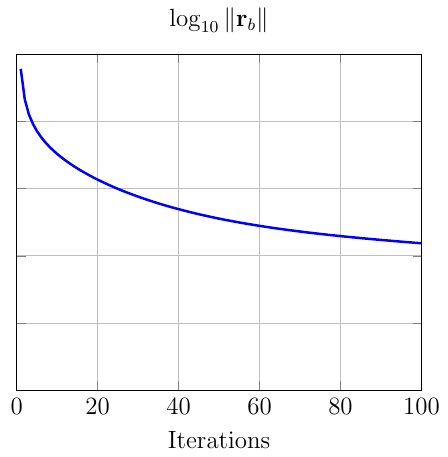}	
\caption{$\log_{10}$ of the Euclidean norm of the residual for the Arnold-Falk-Winther smoother applied to the dual formulation for the Cook's problem in Figure \ref{Cook_problem}. The patch is extended as in Figure \ref{RT1PatchExternalDofs}. Parameters: $N_{\text{dofs}}=4089$, $\mu=1$, $\lambda=1$, $\text{number of smoothing steps}=100$. }
\label{patch_enriched_with_stress_dofs_on_the_boundary_results}
\end{figure}
\begin{figure}[htbp!]
\centering
\includegraphics[width=0.6\textwidth]{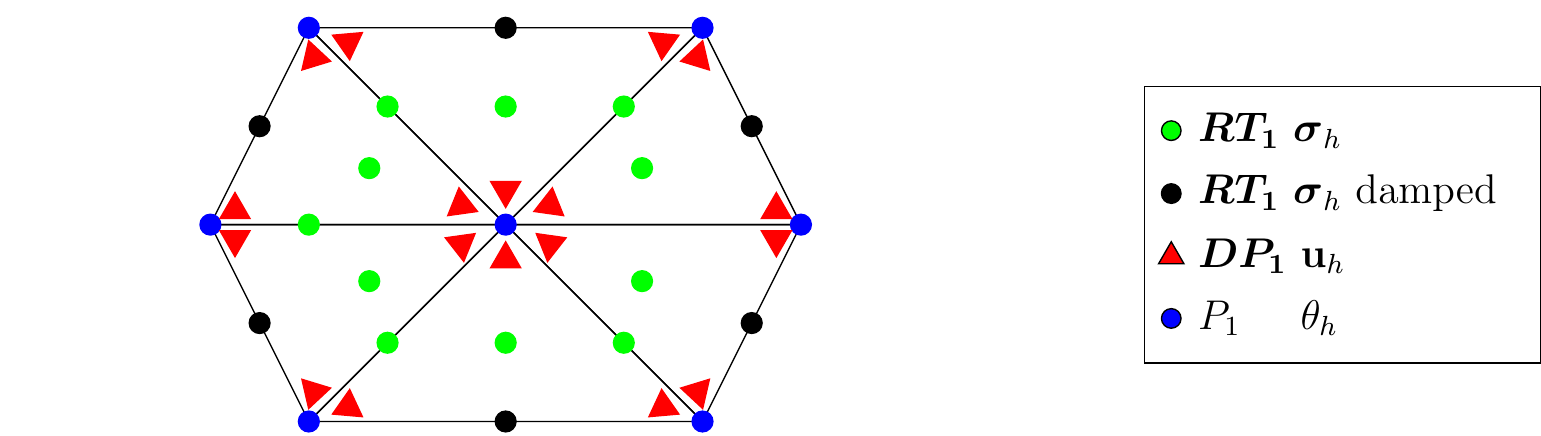}	
\caption{Subspace referred to a patch considering damped $\boldsymbol{RT_1}$ dofs on the border.} 
\label{RT1PatchSampedExternalDofs}
\end{figure}

The subspace of Figure \ref{monolithic_subspace_a} gives rise to a local system that allows for rigid body motions, but that, in principle, satisfies all the constraints belonging to the patch. None of the two proposed strategies, the removal of some dofs or the enforcement of the zero-average constraints, is really able to solve the original local system on the subspace of Figure \ref{monolithic_subspace_a} without important modifications. As a consequence, the corresponding smoothers do not perform well. Another strategy is to enrich the full Neumann problem of Figure \ref{monolithic_subspace_a} by adding also $\bsigma_h$ dofs on the boundary patch, obtaining the full Dirichlet problem of Figure \ref{monolithic_subspace_b}. All the constraints, corresponding to displacement (orange triangles) and rotation (cyan circles) dofs outside the patch but communicating with the $\bsigma_h$ dofs on the border, cannot be fulfilled. See Figure \ref{RT1PatchExternalDofs}. Hence $\norm{\br_b}$ can be expected sometimes to increase after the addition, to the current solution, of the local correction. Nevertheless, since the patches overlap and communicate with each other, it is important the behavior of the residual not after single local corrections, but after a whole smoothing step. In particular, in Figure \ref{patch_enriched_with_stress_dofs_on_the_boundary_results}, we can appreciate that the norm of the residual is decreasing after each smoothing step and, in contrast to Figure \ref{zeroaveragemultipliersresiduals}, the decrease is monotone with a much better rate of decay.\\
The full Dirichlet approach neglects all the dofs of the Lagrange multipliers outside the patch and thus the corresponding constraints cannot be fulfilled. For this reason, it could be profitable to damp the components of the local correction related to the stress boundary dofs: the constraints outside the patch would be violated in a minor way and maybe the convergence could be beneficially affected. The easier way to obtain this penalization is by adding to $\bA_{j,i}$ in (\ref{discrete_ji_saddle_point_system}) a semi-positive definite diagonal matrix that is non-zero only in the positions related to the stress boundary dofs. The new system, built in this way, can be interpreted as an average between the cases of Figure \ref{monolithic_subspace_a} and Figure \ref{monolithic_subspace_b}.\\
The local problem (\ref{discrete_ji_saddle_point_system}) can be rewritten by distinguish between internal and boundary dofs, respectively denoted with the subscripts ``int'' and ``ext''. All other local quantities will be denoted with the subscript ``loc''. The full Dirichlet problem related to Figure \ref{monolithic_subspace_b} becomes:
\begin{align}
\begin{bmatrix}
 \bA_{\text{ext},\text{ext}} &   \bA_{\text{int},\text{ext}}^T & \bB_{\text{int},\text{ext}}^T\\
 \bA_{\text{int},\text{ext}} & \bA_{\text{int},\text{int}} & \bB_{\text{int},\text{int}}^T\\
\bB_{\text{int},\text{ext}} & \bB_{\text{int},\text{int}} & 0
\end{bmatrix}
\begin{bmatrix}
\by_{\text{ext}}\\
\by_{\text{int}}\\
\bz_{\text{loc}}
\end{bmatrix}
=
\begin{bmatrix}
\bff_{\text{ext}}\\
\bff_{\text{int}}\\
\bh_{\text{loc}}
\end{bmatrix}\:,
\label{discrete_local_saddle_point_system_external}
\end{align}
while the full Neumann problem of Figure \ref{monolithic_subspace_a} is:
\begin{align}
\begin{bmatrix}
 \bI&   \textbf{0} & \textbf{0}\\
\textbf{0}& \bA_{\text{int},\text{int}} & \bB_{\text{int},\text{int}}^T\\
\textbf{0 }& \bB_{\text{int},\text{int}} & 0
\end{bmatrix}
\begin{bmatrix}
\by_{\text{ext}}\\
\by_{\text{int}}\\
\bz_{\text{loc}}
\end{bmatrix}
=
\begin{bmatrix}
\textbf{0}\\
\bff_{\text{int}}\\
\bh_{\text{loc}}
\end{bmatrix}\:.
\label{discrete_local_saddle_point_system_internal}
\end{align}
Since the boundary conditions of (\ref{discrete_local_saddle_point_system_internal}) are homogeneous, the problem (\ref{discrete_local_saddle_point_system_internal}) is also equivalent to the following linear system: 
\begin{align}
\begin{bmatrix}
\bG(\alpha)&   \textbf{0} & \textbf{0}\\
\textbf{0}& \bA_{\text{int},\text{int}} & \bB_{\text{int},\text{int}}^T\\
\textbf{0 }& \bB_{\text{int},\text{int}} & 0
\end{bmatrix}
\begin{bmatrix}
\by_{\text{ext}}\\
\by_{\text{int}}\\
\bz_{\text{loc}}
\end{bmatrix}
=
\begin{bmatrix}
\textbf{0}\\
\bff_{\text{int}}\\
\bh_{\text{loc}}
\end{bmatrix}\:,
\label{discrete_local_saddle_point_system_internal_alpha}
\end{align}
where, given the positive scalar $\alpha>0$, $\bG(\alpha)$ is a positive definite diagonal matrix defined in the following way:
\begin{align}
\bG_{p,q}(\alpha)
=\begin{cases}
\alpha  \max \limits_{s}   \max \left(|(\bA_{\text{ext},\text{ext}})_{p,s}| ,|(\bA_{\text{int},\text{ext}}^T)_{p,s}|, |(\bB_{\text{int},\text{ext}}^T)_{p,s}|  \right) & p=q\\
0 & p \neq q
\end{cases}\:.
\label{G_of_alpha_def}
\end{align}
In particular, the system (\ref{discrete_local_saddle_point_system_internal_alpha}) is decoupled: $\by_{\text{ext}}$ and ${[\by_{\text{int}},\bz_{\text{loc}}]^T}$ are independent. For this reason, we can leave the homogeneous boundary conditions as they are and multiply the remaining part of the system by a scalar $\epsilon$. Due to the decoupling, the system is still equivalent for any ${\epsilon>0}$. Thus we can sum up the modified problem (\ref{discrete_local_saddle_point_system_internal_alpha}) and (\ref{discrete_local_saddle_point_system_external}), make $\epsilon$ tend to zero and get the following solvable system:
\begin{align}
\begin{bmatrix}
 \bA_{\text{ext},\text{ext}} +\bG(\alpha)&   \bA_{\text{int},\text{ext}}^T & \bB_{\text{int},\text{ext}}^T\\
 \bA_{\text{int},\text{ext}} & \bA_{\text{int},\text{int}} & \bB_{\text{int},\text{int}}^T\\
\bB_{\text{int},\text{ext}} & \bB_{\text{int},\text{int}} & 0
\end{bmatrix}
\begin{bmatrix}
\by_{\text{ext}}\\
\by_{\text{int}}\\
\bz_{\text{loc}}
\end{bmatrix}
=
\begin{bmatrix}
\bff_{\text{ext}}\\
\bff_{\text{int}}\\
\bh_{\text{loc}}
\end{bmatrix}\:.
\label{discrete_local_saddle_point_system_average}
\end{align}
The new subspace is represented in Figure \ref{RT1PatchSampedExternalDofs}.

Intuitively, in the new problem (\ref{discrete_local_saddle_point_system_average}), discrete Robin conditions are enforced, because (\ref{discrete_local_saddle_point_system_average}) is recovered as the average of full Neumann and full Dirichlet problems. The choice of using Robin, instead of Dirichlet or Neumann boundary conditions, in a domain decomposition approach for primal formulations, has been proved to be an efficient choice. See \cite{EG03}, \cite{Gan08}, \cite{GHM07}, \cite{GHM+07}, \cite{GAT07}, \cite{GV19}. \\
In order to determine the relation between (\ref{discrete_local_saddle_point_system_average}) and the case of discretized Robin conditions, let us consider, for the sake of simplicity,  the strong formulation of the Poisson problem:
\begin{subequations}
\begin{equation}
\begin{alignedat}{2}
-\text{div}\bsigma &= f &&\quad \text{on } \Omega \:,\\
\bsigma&=\nabla u &&\quad \text{on } \Omega \:,
\end{alignedat}
\end{equation} 
\end{subequations}
where $f \in L^2(\Omega)$. We multiply the equilibrium equation by a test function $v \in L^2(\Omega)$ and the constitutive equation by a test $\btau \in  H_{\text{div}}(\Omega)$. Then we integrate the second equation by parts and obtain:
\begin{subequations}
\begin{equation}
\begin{alignedat}{2}
\int_{\Omega} \bsigma \cdot \btau +\int_{\Omega} u  \: \text{div} \btau-\int_{\partial \Omega}u\: \left(\btau \cdot \bn \right) &=0   &&\quad \forall \btau \in  H_{\text{div}}(\Omega)\:, \\
\int_{\Omega}  \text{div} \bsigma \: v&= \int_{\Omega} f v  &&\quad \forall v \in L^2(\Omega)\:.
\end{alignedat}
\end{equation} 
\end{subequations}
By enforcing homogeneous Robin boundary conditions for the dual formulation, i.e:
\begin{align}
u+\alpha \bsigma \cdot \bn =0\:,
\end{align}
we seek for $(u,\bsigma) \in L^2(\Omega) \times  H_{\text{div}}(\Omega)$ such that:
\begin{subequations}
\begin{equation}
\begin{alignedat}{2}
\int_{\Omega} \bsigma \cdot \btau +\int_{\Omega} u  \: \text{div} \btau +\int_{\partial \Omega} \alpha\:   \left(\bsigma \cdot \bn \right)\: \left(\btau \cdot \bn \right) &=0   &&\quad \forall \btau \in  H_{\text{div}}(\Omega)\:, \\
\int_{\Omega}  \text{div} \bsigma \: v&= \int_{\Omega} f v  &&\quad \forall v \in L^2(\Omega)\:.
\end{alignedat}
\end{equation} 
\end{subequations}
The problem is discretized by means of Raviart-Thomas and discontinuous Lagrangian finite elements so that an LBB condition is satisfied. All the boundary and all the internal dofs of $\bsigma$ are collected in $\bsigma_{\text{ext}}$ and in $\bsigma_{\text{int}}$, while  all the dofs of $\bu$ are inserted in $\bu_{\text{loc}}$. The following system is obtained:
\begin{align}
\begin{bmatrix}
\bS_{\text{ext},\text{ext}}+\alpha \bM & \bS_{\text{int},\text{ext}}^T & \bT_{\text{ext}}^T\\
\bS_{\text{int},\text{ext}} & \bS_{\text{int},\text{int}} & \bT_{\text{int}}^T\\
\bT_{\text{ext}} & \bT_{\text{int}}       & \textbf{0}
\end{bmatrix}
\begin{bmatrix}
\bsigma_{\text{ext}}\\
\bsigma_{\text{int}}\\
\bu_{\text{loc}}
\end{bmatrix}
=
\begin{bmatrix}
\textbf{0}\\
\textbf{0}\\
\bff_{\text{loc}}
\end{bmatrix}
\label{full_robin_poisson_discrete_dual}
\end{align}
where the matrices $\bS$, $\bM$ and $\bT$ are respectively the discretizations of the bilinear forms ${s(\bsigma,\btau)=\int_{\Omega} \bsigma \cdot \btau}$, ${m(\bsigma,\btau)=\int_{\partial \Omega} \alpha\:   \left(\bsigma \cdot \bn \right)\: \left(\btau \cdot \bn \right) }$ and $t(\bsigma,v)=\int_{\Omega}\text{div} \bsigma \: v$, while the vector $\bff_{\text{loc}}$ discretizes $f(v)=\int_{\Omega} f\:v$. As in the primal case, also in the dual formulation Robin conditions give rise to a mass matrix on the border, $\bM$. In principle, some components of the volumetric mass matrix obtained by means of Raviart-Thomas functions can be negative. However, the trace of the Raviart-Thomas shape functions for simplicial meshes reduces to discontinuous Lagrange finite elements on the border. For this reason, the mass matrix $\bM$ on the border is still similar to the one of the primal case.\\
In order to recover the conditions in (\ref{discrete_local_saddle_point_system_average}) for the dual Poisson problem, the two discrete problems, with homogenous Dirichlet and homogenous Neumann boundary conditions, have to be introduced. For the full Dirichlet case, we get:
\begin{align}
\begin{bmatrix}
\bS_{\text{ext},\text{ext}} & \bS_{\text{int},\text{ext}}^T & \bT_{\text{ext}}^T\\
\bS_{\text{int},\text{ext}} & \bS_{\text{int},\text{int}} & \bT_{\text{int}}^T\\
\bT_{\text{ext}} & \bT_{\text{int}}       & 0
\end{bmatrix}
\begin{bmatrix}
\bsigma_{\text{ext}}\\
\bsigma_{\text{int}}\\
\bu_{\text{loc}}
\end{bmatrix}
=
\begin{bmatrix}
\textbf{0}\\
\textbf{0}\\
\bff_{\text{loc}}
\end{bmatrix}\:,
\label{full_dirichlet_poisson_discrete_dual}
\end{align}
while for Neumann boundary conditions, we get:
\begin{align}
\begin{bmatrix}
\bG & \textbf{0} & \textbf{0}\\
\textbf{0} & \bS_{\text{int},\text{int}} & \bT_{\text{int}}^T\\
\textbf{0}& \bT_{\text{int}}       & 0
\end{bmatrix}
\begin{bmatrix}
\bsigma_{\text{ext}}\\
\bsigma_{\text{int}}\\
\bu_{\text{loc}}
\end{bmatrix}
=
\begin{bmatrix}
\textbf{0}\\
\textbf{0}\\
\bff_{\text{loc}}
\end{bmatrix}\:.
\label{full_neumann_poisson_discrete_dual}
\end{align}
With the same trick explained used in \ref{discrete_local_saddle_point_system_internal_alpha}, we average between 
(\ref{full_dirichlet_poisson_discrete_dual})
and (\ref{full_neumann_poisson_discrete_dual}), obtaining:
\begin{align}
\begin{bmatrix}
\bS_{\text{ext},\text{ext}}+ \bG & \bS_{\text{int},\text{ext}}^T & \bT_{\text{ext}}^T\\
\bS_{\text{int},\text{ext}} & \bS_{\text{int},\text{int}} & \bT_{\text{int}}^T\\
\bT_{\text{ext}} & \bT_{\text{int}}      & 0
\end{bmatrix}
\begin{bmatrix}
\bsigma_{\text{ext}}\\
\bsigma_{\text{int}}\\
\bu_{\text{loc}}
\end{bmatrix}
=
\begin{bmatrix}
\textbf{0}\\
\textbf{0}\\
\bff_{\text{loc}}
\end{bmatrix}\:.
\label{full_robin_poisson_average_discrete_dual}
\end{align}
By defining each diagonal entry of $\bG$ as a scaled sum of the corresponding row in $\bM$:
\begin{align}
\bG_{i i}= \left(\dfrac{ \bG_{i i}} {\sum_{j} \bM_{i j}} \right)\sum_{j} \bM_{i,j} = \beta_i \sum_{j} \bM_{i,j}\:,
\end{align}
the matrix $\bG$ now represents a lumping of $\bM$, scaled by the coefficients $\beta_i$. Therefore the discrete Robin conditions in (\ref{full_robin_poisson_average_discrete_dual}) represent the lumped discretized version of the continuous Robin conditions in (\ref{full_robin_poisson_discrete_dual}), where varying coefficients $\beta_i$ are used instead of $\alpha$. The algebraic representation of (\ref{full_robin_poisson_average_discrete_dual}) is nonethless much easier to implement for different fine and coarse patches than its exact formulation (\ref{full_robin_poisson_discrete_dual}).

\section{Numerical examples in dual linear elasticity}
\label{numerical_examples_linear_dual_elasticity_section}
In this section, the convergence behavior of the patch smoother with full Robin conditions is examined. Due to the equation (\ref{discrete_local_saddle_point_system_average}), the smoother now depends on an input parameter $\alpha \geq 0$. For $\alpha=0$, we recover the full Dirichlet case. For $\alpha > 0$, the smoother will behave differently and here we want to examine how it can affect also the multigrid method convergence.
\subsection{Multigrid for the Cook's membrane problem}
\begin{figure}[htbp!]
\hspace{+3cm}
\includegraphics[width=0.267\textwidth]{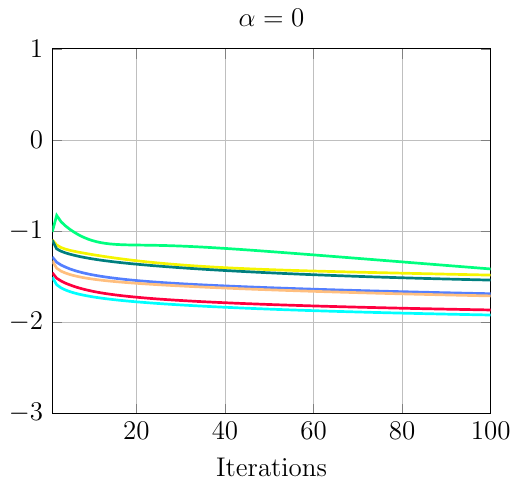}	
\includegraphics[width=0.25\textwidth]{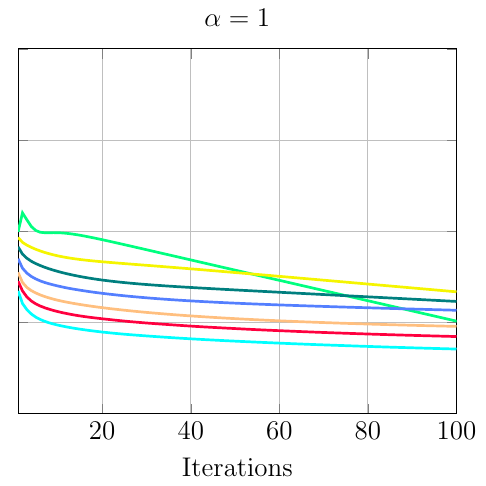}

\hspace{+3cm}\includegraphics[width=0.267\textwidth]{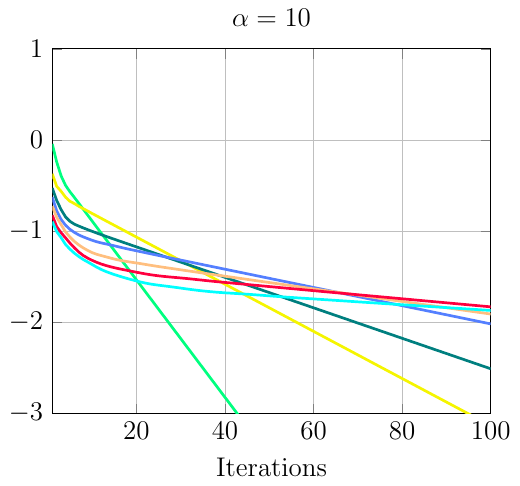}		
\includegraphics[width=0.25\textwidth]{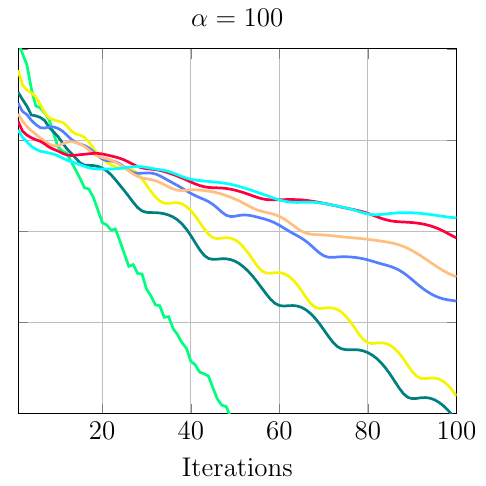}
\raisebox{3.5ex}{
\includegraphics[width=0.15\textwidth]{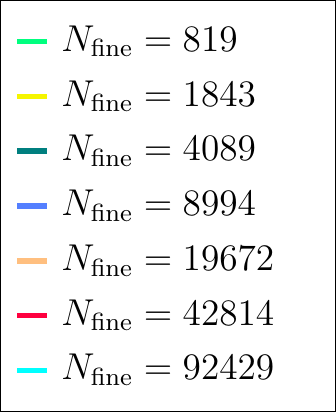}	}	
\caption{$\log_{10}$ of the Euclidean norm of the residual for the patch monolitich smoother applied to the dual formulation for the problem in Figure \ref{Cook_problem}. Parameters: $\mu=1$, $\lambda=\infty$. The residuals have been computed after each smoothing step. }
\label{Cook_smoothing_alpha}
\end{figure}

\begin{figure}[htbp!]
\hspace{+3cm}
\includegraphics[width=0.267\textwidth]{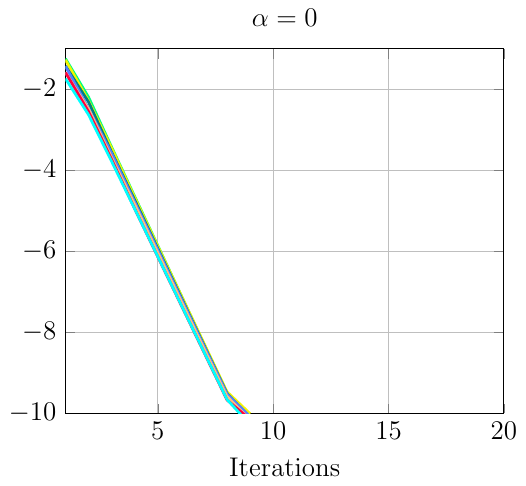}	
\includegraphics[width=0.243\textwidth]{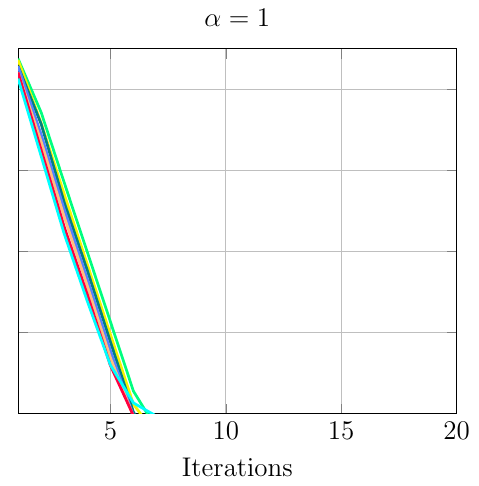}

\hspace{+3cm}
\includegraphics[width=0.267\textwidth]{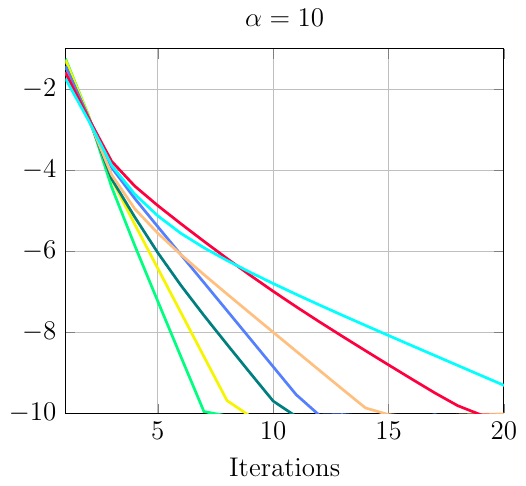}		
\includegraphics[width=0.243\textwidth]{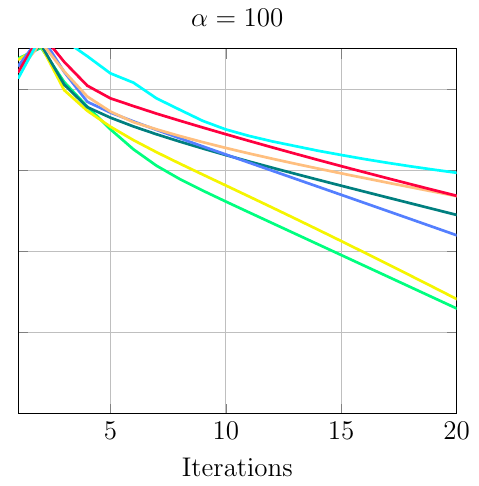}
\raisebox{2.5ex}{
\includegraphics[width=0.15\textwidth]{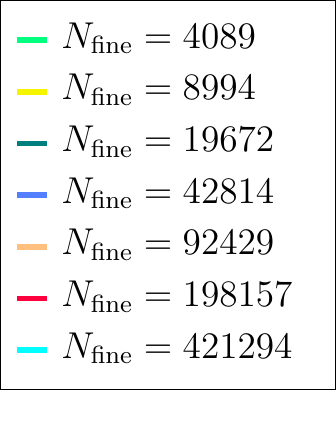}	}
\caption{$\log_{10}$ of the Euclidean norm of the residual for the multigrid method applied to the dual formulation for the Cook's problem in Figure \ref{Cook_problem}. Parameters: $\mu=1$, $\lambda=\infty$, $\text{number of smoothing steps}=5$. The coarsest level has dimension $N_{\text{coarse}}=82$. Then bisection on each element is used to refine the mesh. The residuals have been computed after each V-cycle. }
\label{Cook_multigrid_alpha}
\end{figure}

\begin{figure}[htbp!]
\hspace{+3cm}
\includegraphics[width=0.267\textwidth]{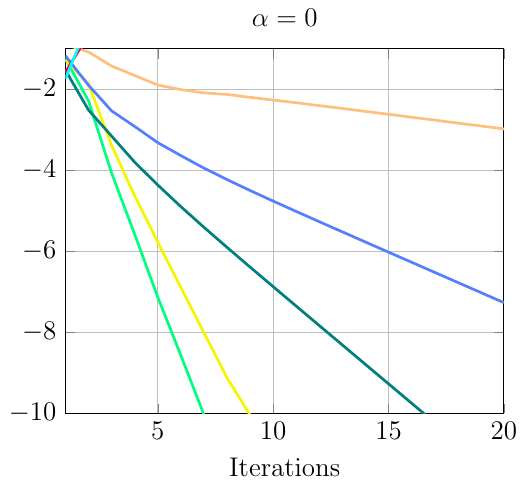}	
\includegraphics[width=0.243\textwidth]{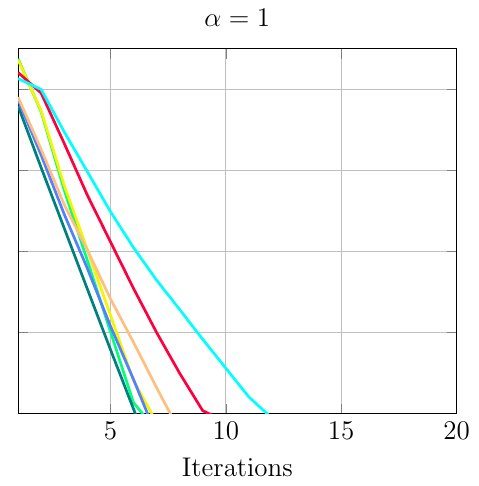}

\hspace{+3cm}\includegraphics[width=0.267\textwidth]{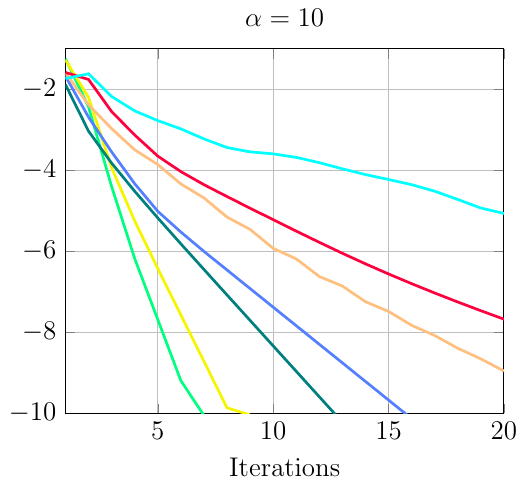}		
\includegraphics[width=0.243\textwidth]{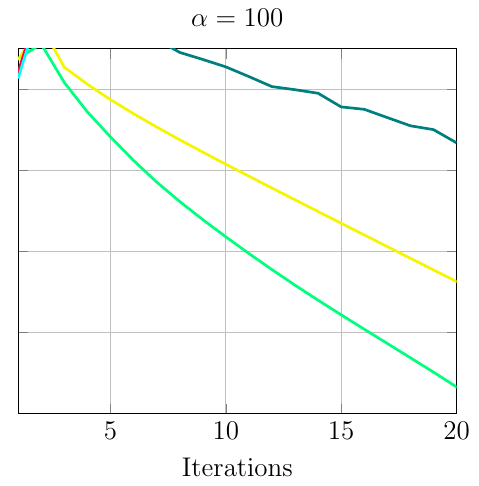}
\raisebox{2.3ex}{
\includegraphics[width=0.15\textwidth]{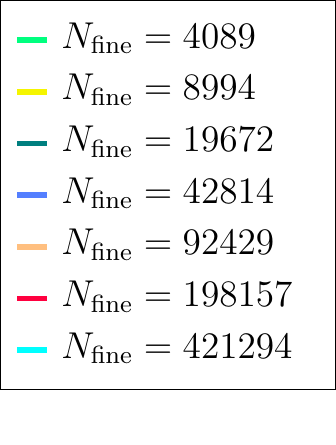}	}
\caption{$\log_{10}$ of the Euclidean norm of the residual for the 2 grids method applied to the dual formulation for the Cook's problem in Figure \ref{Cook_problem}. Parameters: $\mu=1$, $\lambda=\infty$, $\text{number of smoothing steps}=5$. The coarsest level has dimension $N_{\text{coarse}}=819$. Then bisection on each element is used to refine the mesh, but only the finest level is considered. The residuals have been computed after each V-cycle. }
\label{Cook_2grid_alpha}
\end{figure}

\begin{figure}[htbp!]
\centering
\includegraphics[width=0.215\textwidth]{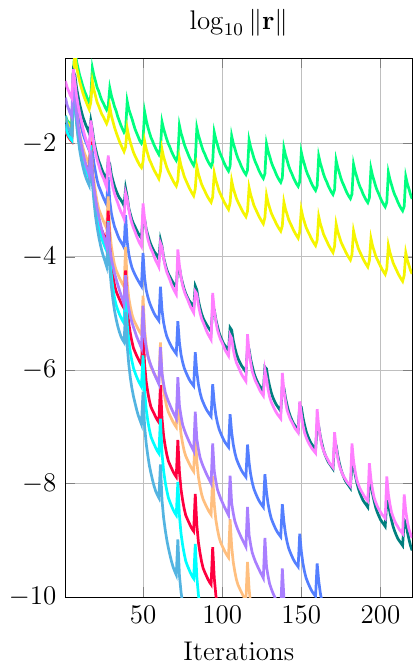}	
\includegraphics[width=0.19\textwidth]{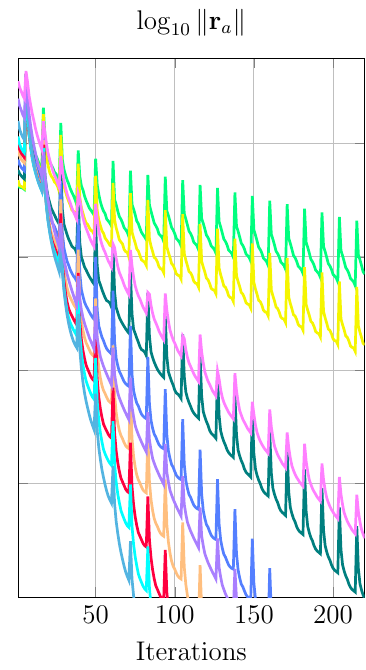}
\includegraphics[width=0.19\textwidth]{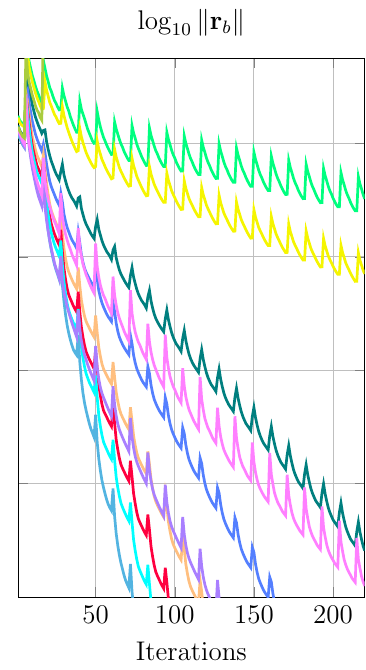}	
\quad
\raisebox{2.ex}{
\includegraphics[width=0.1\textwidth]{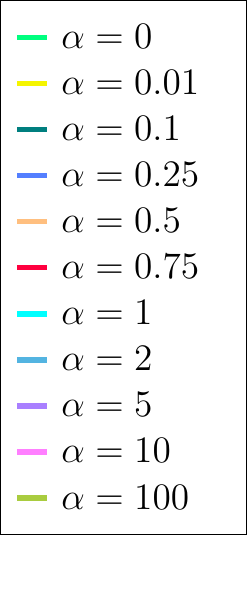}	}
\caption{
$\log_{10}$ of the Euclidean norm of the residual for the 2 grids method applied to the dual formulation for the Cook's problem in Figure \ref{Cook_problem}. Parameters: $N_{\text{coarse}}=819$, $N_{\text{fine}}=92429$, $\mu=1$, $\lambda=\infty$, $\text{number of smoothing steps}=5$. The residuals have been computed after each smoothing step and each coarse correction addition to the current solution. In particular, after the coarse correction addition, the residual suddenly increases. This is due to the coarse representation of equality constraints. Nevertheless the coarse correction still ensures a faster convergence than for the only smoother case of Figure \ref{Cook_smoothing_alpha}.
}
\label{Cook_multigrid_alpha_oscillations}
\end{figure}

\begin{figure}[htbp!]

\hspace{+2cm}
\begin{subfigure}{.35\textwidth}
\includegraphics[scale=0.2]{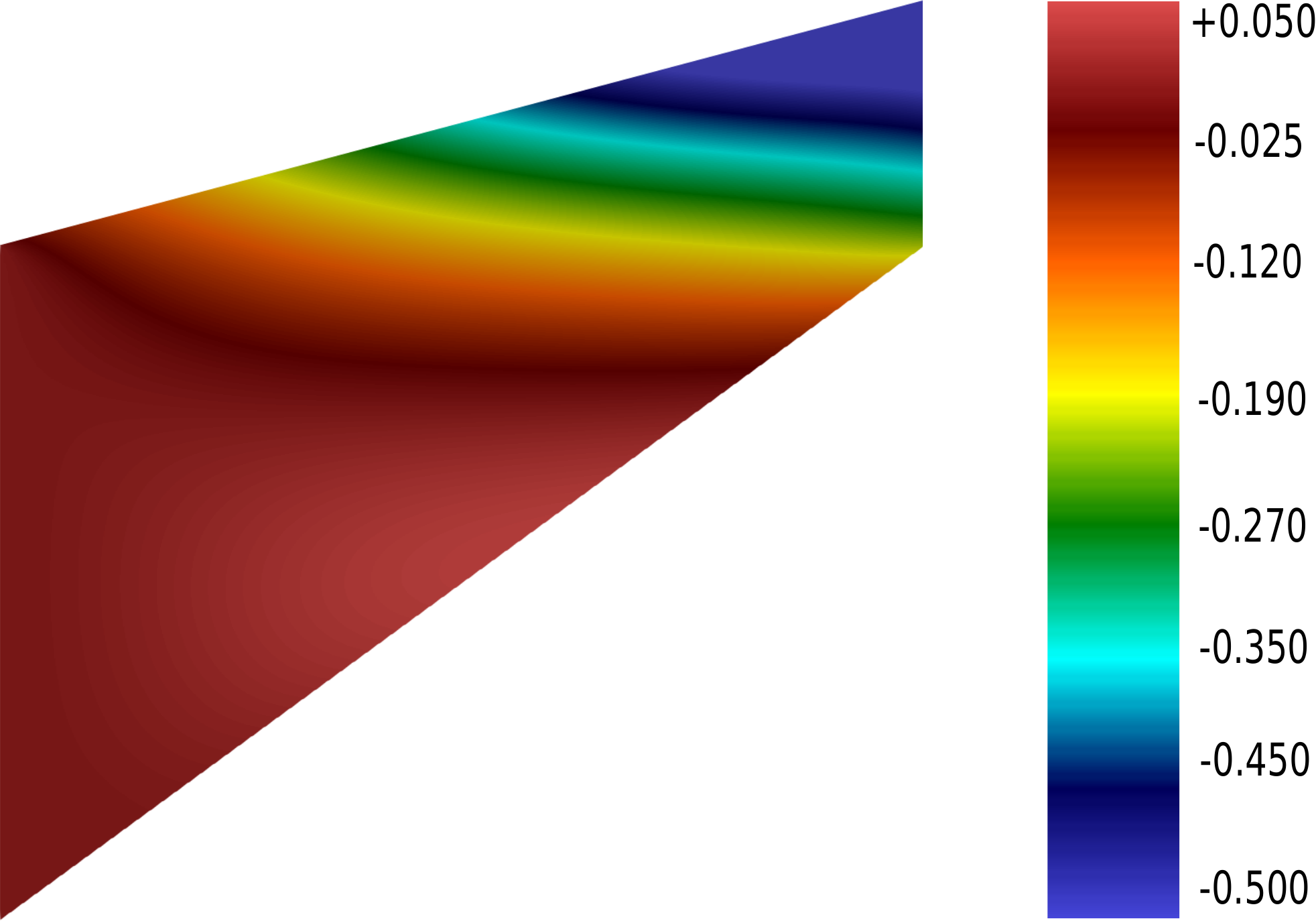}
\caption{$u_x\qquad \qquad\qquad\qquad\qquad\quad$                                 }
\end{subfigure}
$\qquad$ $\qquad$
\begin{subfigure}{.31\textwidth}
\includegraphics[scale=0.2]{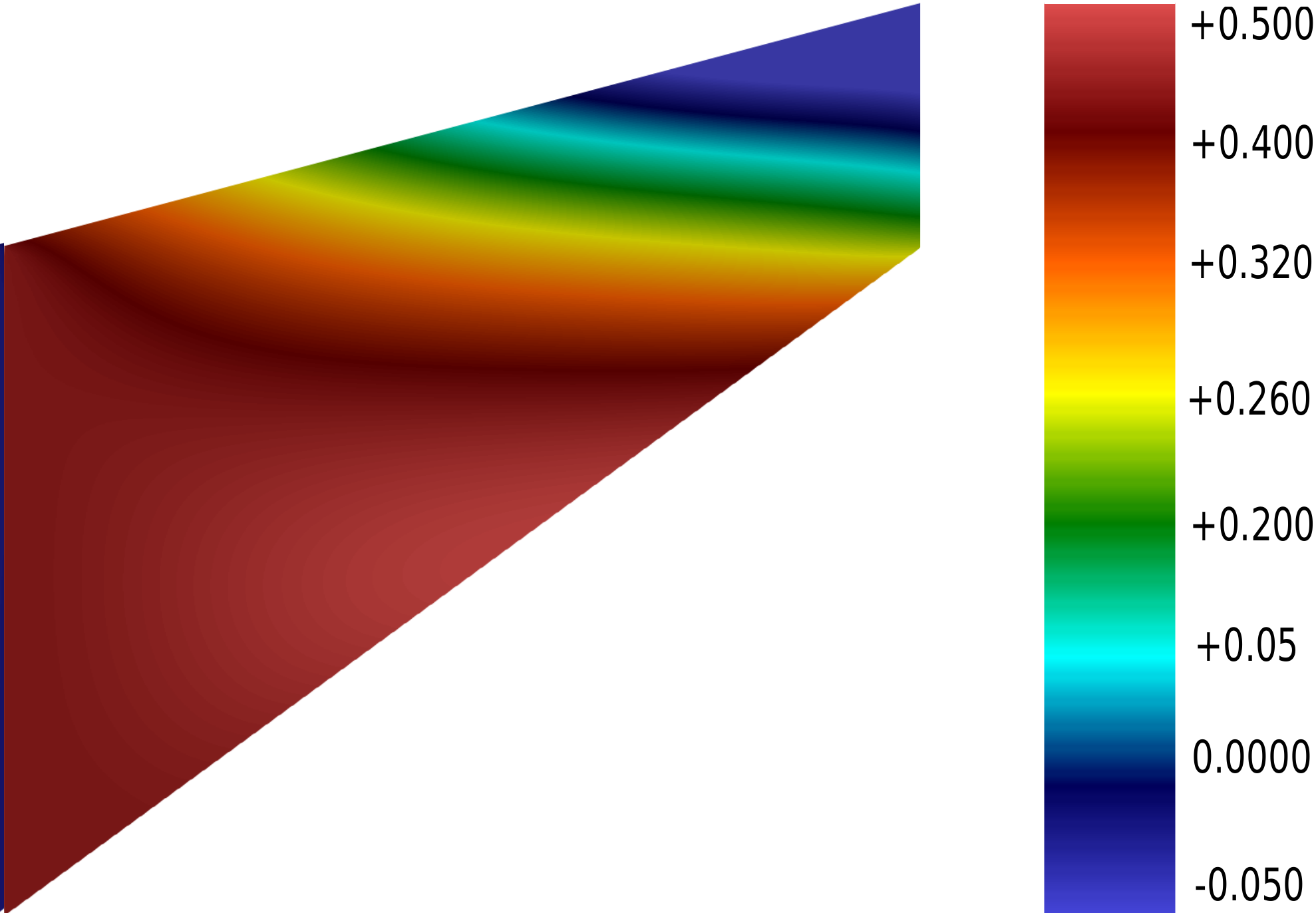}
\caption{$u_y\qquad \qquad\qquad\qquad\qquad$                                 }
\end{subfigure}

\hspace{+2cm}
\begin{subfigure}{.35\textwidth}
\includegraphics[scale=0.2]{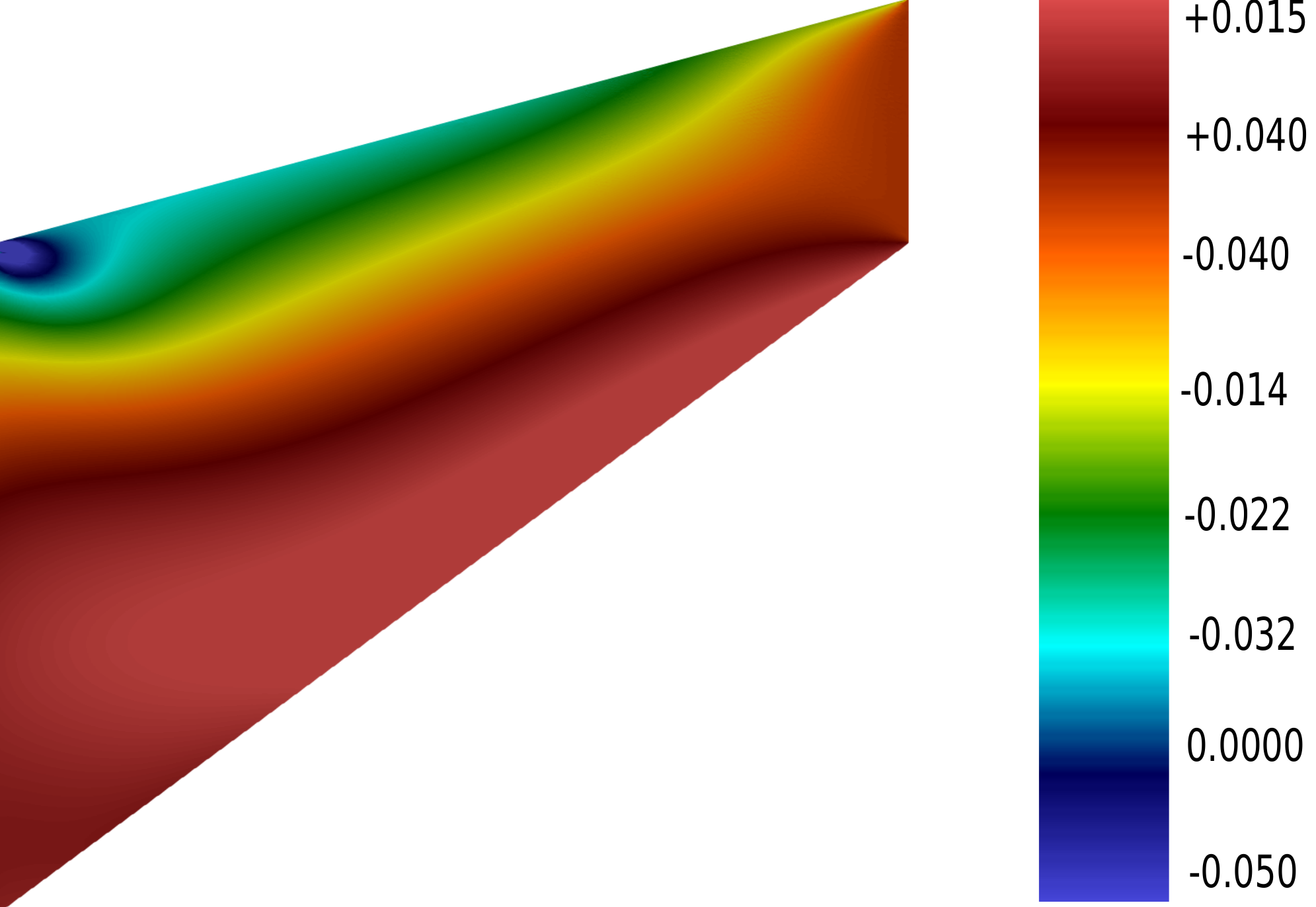}
\caption{$\sigma_{xx}\qquad \qquad\qquad\qquad\qquad\quad$                                             }
\end{subfigure}
$\qquad$ $\qquad$
\begin{subfigure}{.35\textwidth}
\includegraphics[scale=0.2]{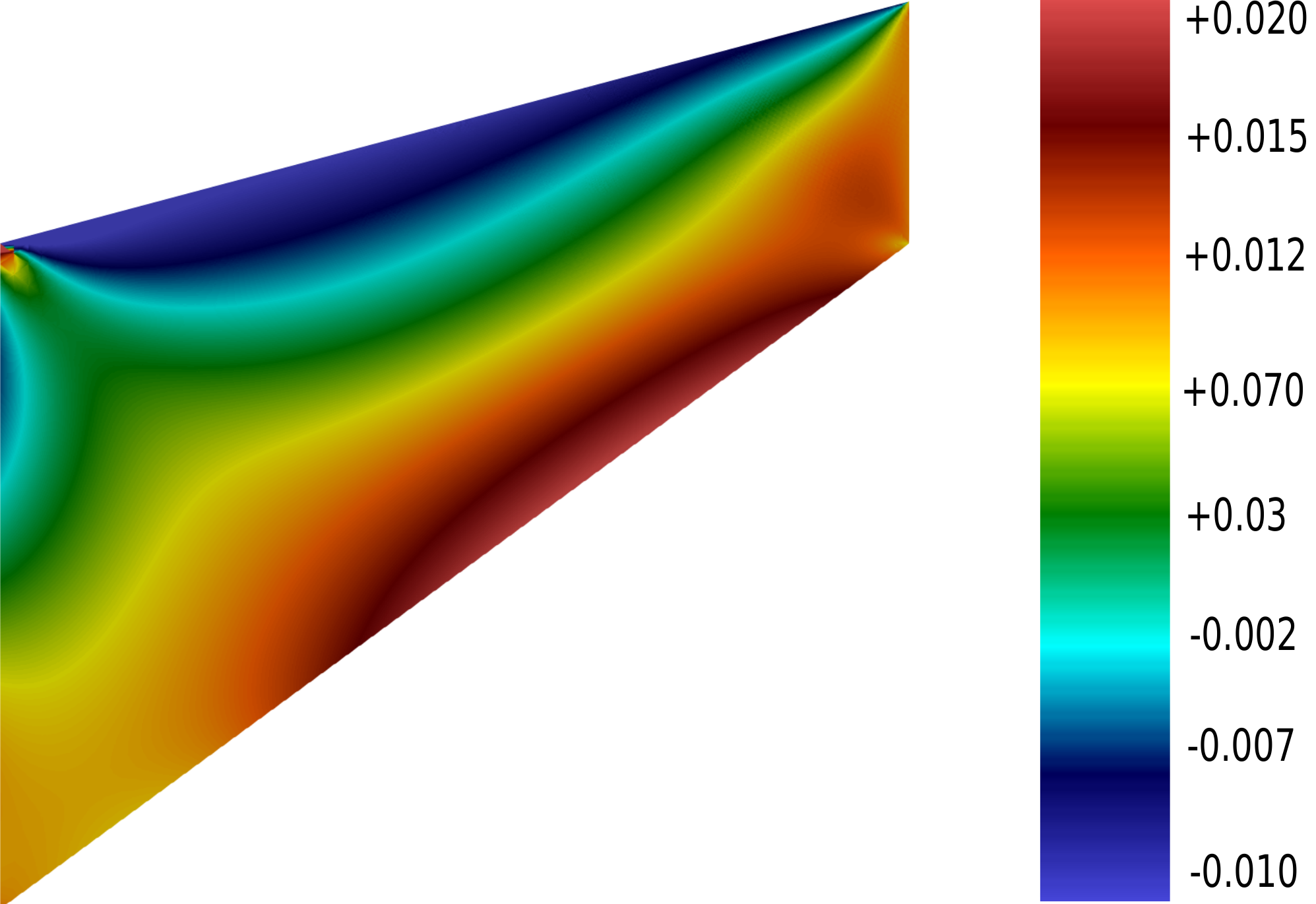}
\caption{$\sigma_{xy}\qquad \qquad\qquad\qquad\qquad \quad$                                 }
\end{subfigure}

\hspace{+2cm}
\begin{subfigure}{.35\textwidth}
\includegraphics[scale=0.2]{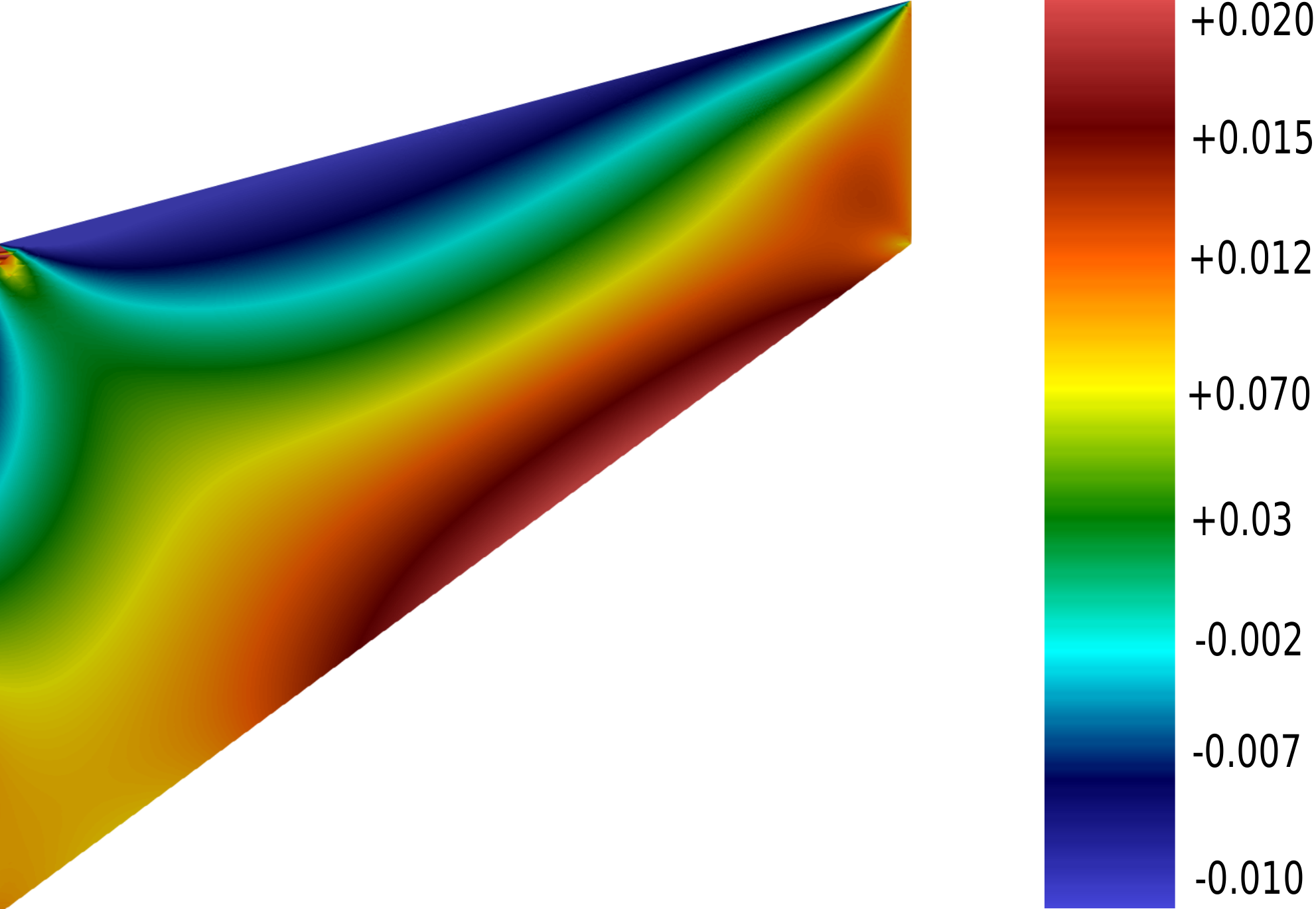}
\caption{$\sigma_{yx}\qquad \qquad\qquad\qquad\qquad\quad$                                 }
\end{subfigure}
$\qquad$ $\qquad$
\begin{subfigure}{.31\textwidth}
\includegraphics[scale=0.2]{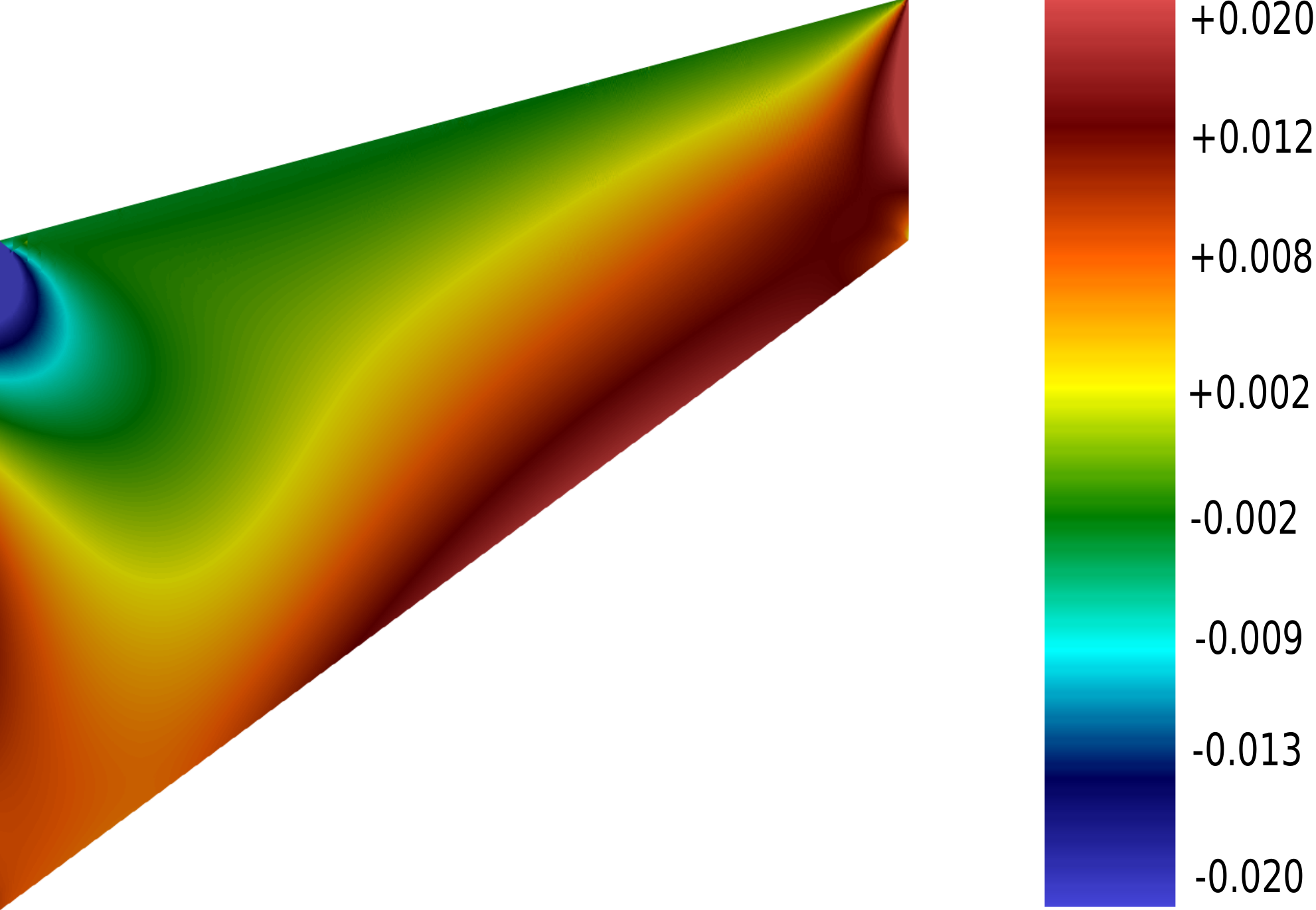}
\caption{$\sigma_{yy}\qquad \qquad\qquad\qquad\qquad$                                            }
\end{subfigure}

\hspace{+2cm}
\begin{subfigure}{.35\textwidth}
\includegraphics[scale=0.2]{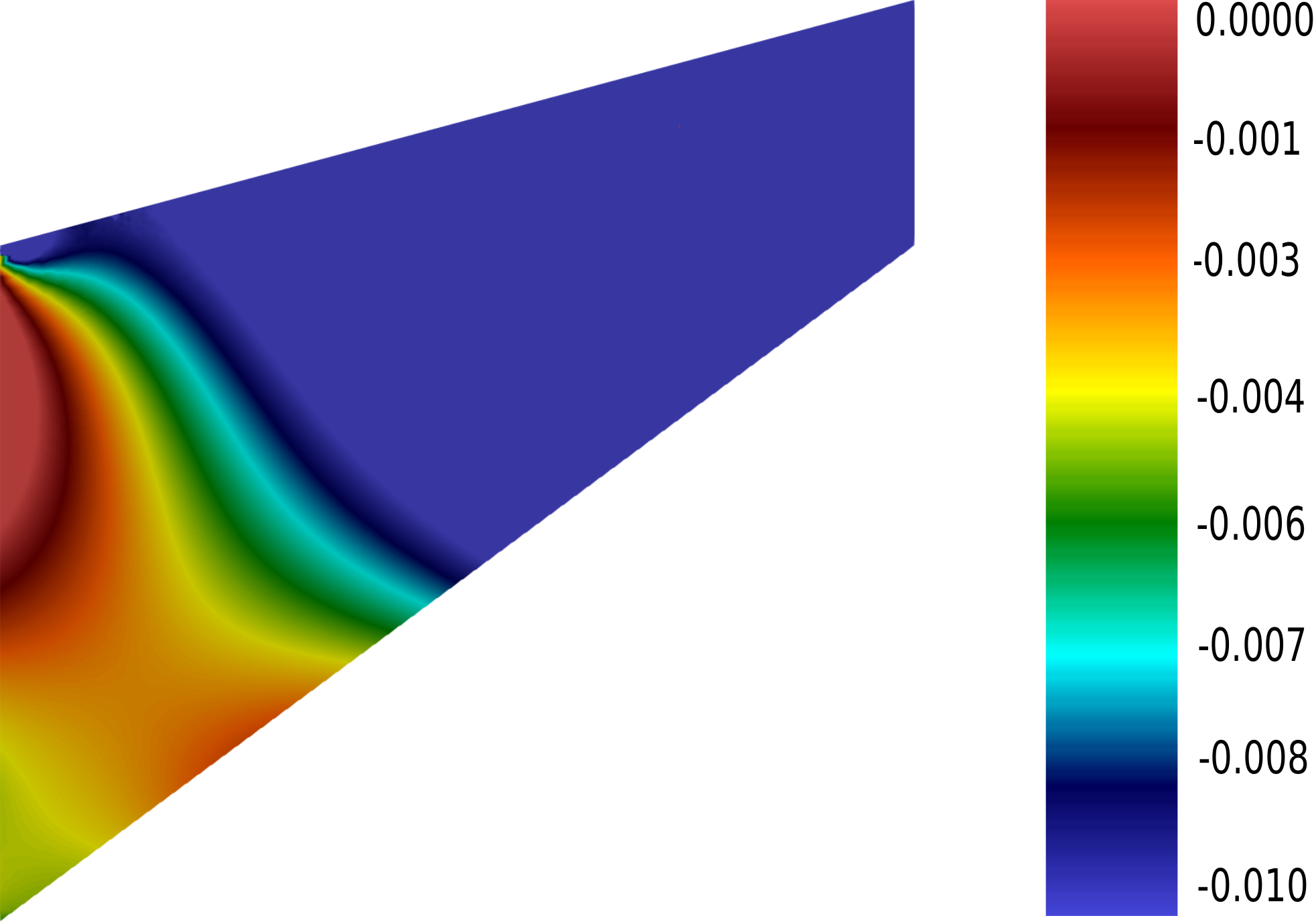}
\caption{$\rho\qquad \qquad\qquad\qquad\qquad\quad$                                              }
\end{subfigure}
$\qquad$ $\qquad$
\begin{subfigure}{.35\textwidth}
\includegraphics[scale=0.075]{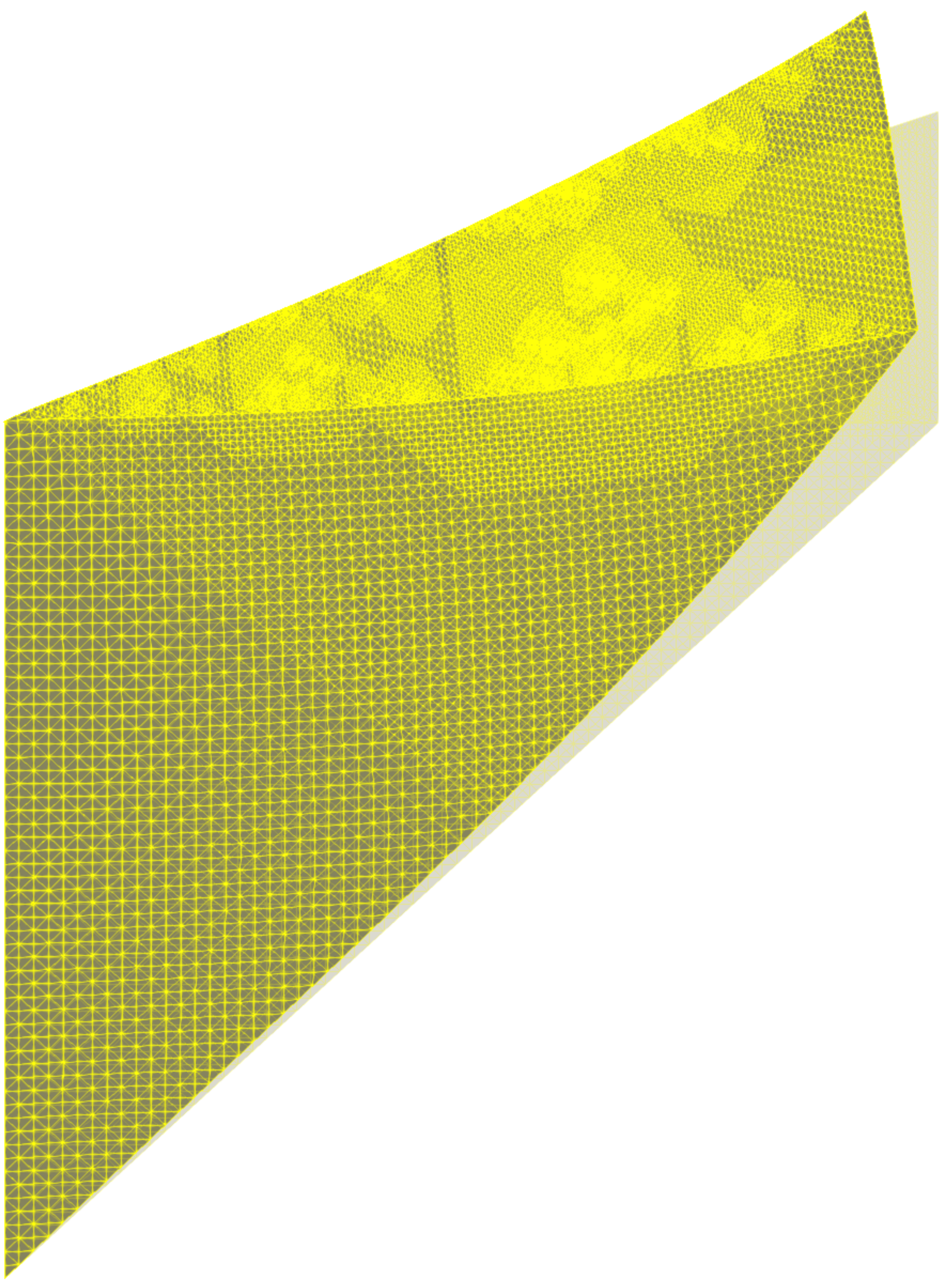}
\caption{5 times larger deformed mesh $\qquad \qquad $ }
\end{subfigure}
\caption{Results for the Cook's problem in Figure \ref{Cook_problem}. Parameters: $\mu=1$, $\lambda=\infty$.}
\label{cook_solutions}
\end{figure}

We now consider the Cook's problem, as depicted in Figure \ref{Cook_problem}. On the left edge, we enforce a zero displacament. On the right edge, a vertical force is applied: $\bsigma \cdot \bn=[0,0.01]^T$. Everywhere else, we impose homogeneous Neumann conditions, $\bsigma \cdot \bn =\textbf{0}$. The material has the following parameters: $\mu=1$ and $\lambda=\infty$. We have analyzed the cases for different values of $\alpha$, but the more significative ones for this experiment are ${\alpha=0, \: 1, \: 10, \: 100}$. The solutions can be found in Figure \ref{cook_solutions}. \\
Figure \ref{Cook_smoothing_alpha} illustrates a convergent behavior of the smoother for different values of $\alpha$. The parameter ${0 \leq \alpha \leq 1}$ does not seem to affect too much the performance. For ${\alpha=10,\:100}$, the convergence is faster for small meshes, but it gets worse for larger meshes. So by inspecting the only smoother, we can assume the multigrid method will not be affected for ${\alpha \in [0,1]}$. \\
We use the multigrid method with a coarse mesh of $N_{\text{coarse}}=819$ dofs which we refine, with a bisection algorithm, up to $N_{\text{fine}}=421294$ dofs. For each level, we do 5 pre-smoothing steps and 5 post-smoothing steps. On the coarsest level, we solve exactly. Again, we have repeated the experiments for the same values of  $\alpha$. Similarly to Figure \ref{Cook_smoothing_alpha}, Figure \ref{Cook_multigrid_alpha} shows a rate of convergence that is not much influenced by $\alpha \in [0,1]$. In this case, optimal convergence is achieved, meaning that the number of iterations is independent of the number of dofs and the number of levels used. On the other hand, for $\alpha>1$, the rate of convergence is no more optimal and depends on the size of the problem. For too large values of $\alpha$, the method does not even converge.\\
The impact of $\alpha > 0$ does not seem so necessary for a linear multigrid method that considers all the levels. However, if we inspect a 2-grid method as in Figure \ref{Cook_2grid_alpha}, its role becomes more and more important with the aggressivity of the coarsening. The best performance is attained for $\alpha=1$. We recover a convergence rate that is almost independent of the dimension of the problem, even though only two levels are used. Thus, if the right value of $\alpha$ is chosen, the Robin boundary conditions can damp different frequency components of the error. As we move away from this value, the method becomes slower and slower and can also not converge. The choice $\alpha=0$ for too fine meshes does not make the 2-grids method convergent. We finally understand that the choice of $\alpha$ is very important and could affect more complicated non-linear problems.\\
In Figure \ref{Cook_multigrid_alpha_oscillations}, we see how the norm of the residual behaves after each fine and coarse correction for a 2-grid method. This means we evaluate the global fine residual after the smoothing steps and after the coarse correction. Every time the latter is added to the current solution, the norm of the residual increases. Indeed the fine equality constraints are just projected onto the coarser space. Then on the coarse space, an exact coarse correction is computed and it is interpolated to the fine space. Such correction, however, just satisfies the coarse and not the fine constraints. This is why, after its addition to the current solution, we can see an increment of the norm of the residual. Nevertheless, at the same time, the coarse correction helps in the global communication process and, thanks to the post-smoothing steps, in accelerating the overall convergence. In conclusion, we can state that the parameter $\alpha$ does not only govern the communication among the subdomains, but it is also important for damping the error after the computation of coarse corrections that do not fully satisfy the fine constraints. In particular, the value of $\alpha$ is properly chosen for the convergence if $\norm{\br_a}$ and $\norm{\br_b}$ are comparable and the latter does not dominate the whole process.

\clearpage
\subsection{Multigrid for the face problem}
\begin{figure}[htbp!]
\centering
\includegraphics[width=0.267\textwidth]{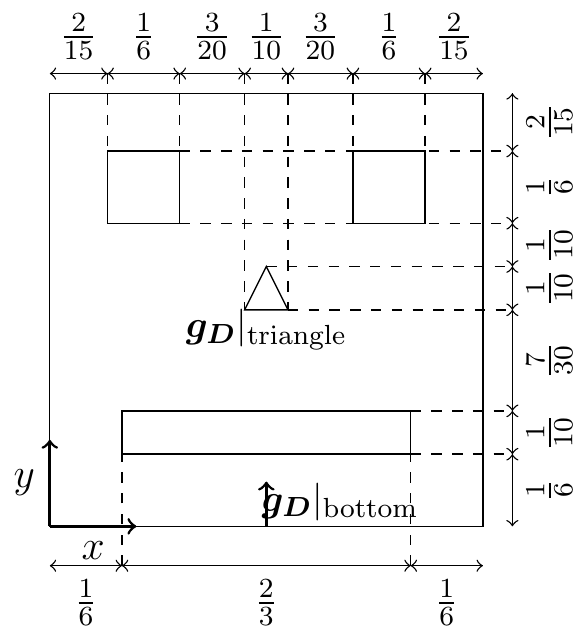}	
\caption{The face geometry is a square domain ${[0,1] \times [0,1]}$ with four holes, three rectangles and one triangle. The displacament condition on the bottom side is ${\boldsymbol{g_{D}}|_{\text{bottom}}=[0,0.05 x^2]^T}$. On the triangle-shaped hole ${\boldsymbol{g_D}|_{\text{triangle}}=[0,0]^T}$. Everywhere else ${\boldsymbol{g_N}=[0,0]^T}$.}
\label{face_problem}
\end{figure}

\begin{figure}[htbp!]
\hspace{+3cm}
\includegraphics[width=0.27\textwidth]{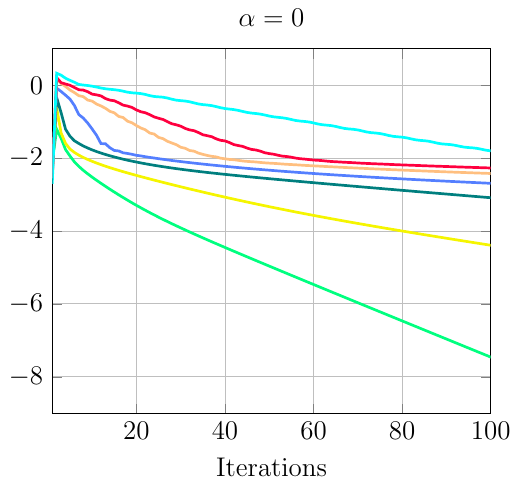}
\includegraphics[width=0.25\textwidth]{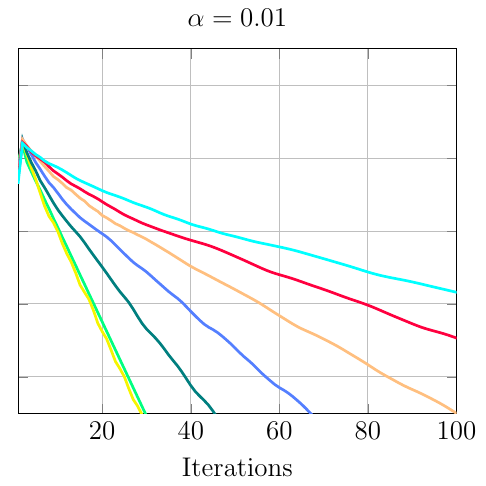}

\hspace{+3cm}\includegraphics[width=0.27\textwidth]{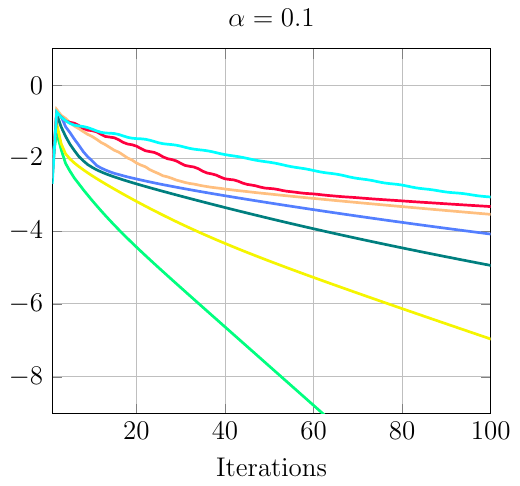}
\includegraphics[width=0.25\textwidth]{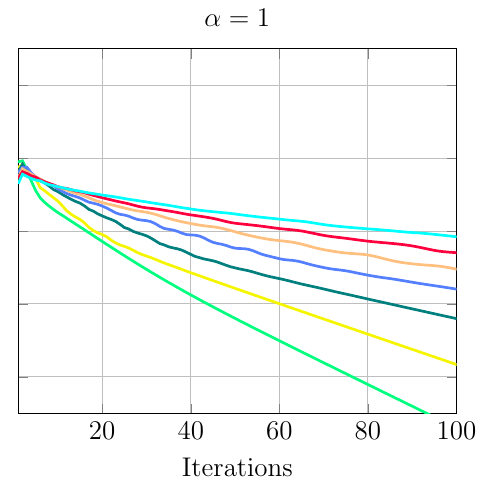}	
\raisebox{3.5ex}{
\includegraphics[width=0.15\textwidth]{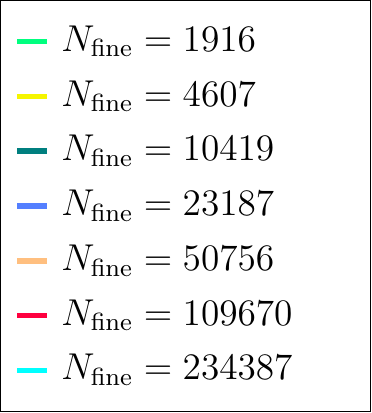}	}
\caption{$\log_{10}$ of the Euclidean norm of the residual for the the patch monolitich smoother applied to the dual formulation for the problem in Figure \ref{face_problem}. Parameters: $\mu=1$, $\lambda=\infty$. The residuals have been computed after each smoothing step. }
\label{face_smoothing_alpha}
\end{figure}

\begin{figure}[htbp!]
\hspace{+3cm}
\includegraphics[width=0.27\textwidth]{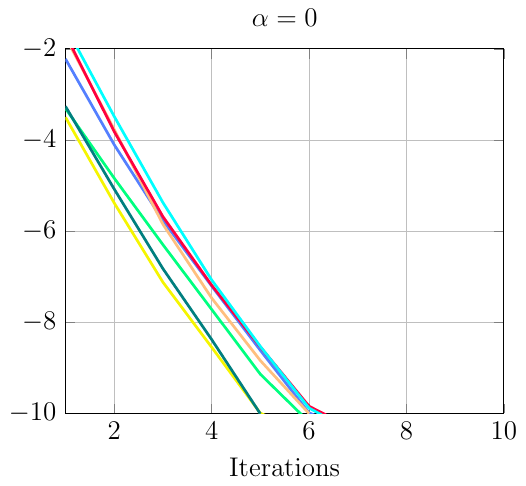}
\includegraphics[width=0.25\textwidth]{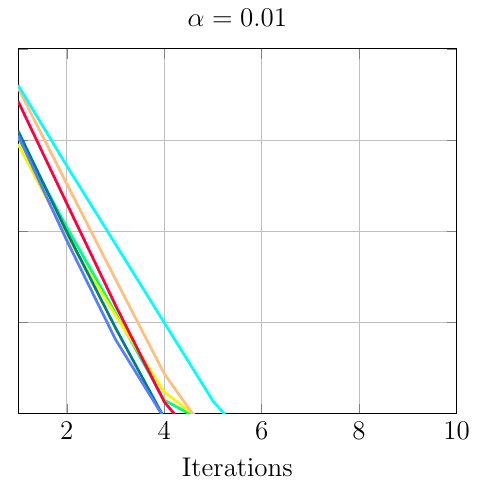}

\hspace{+3cm}\includegraphics[width=0.27\textwidth]{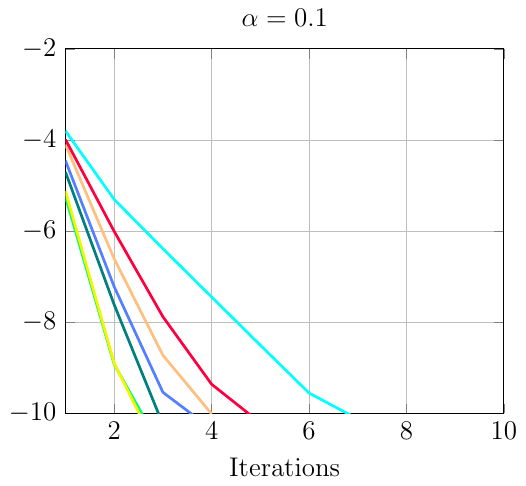}
\includegraphics[width=0.25\textwidth]{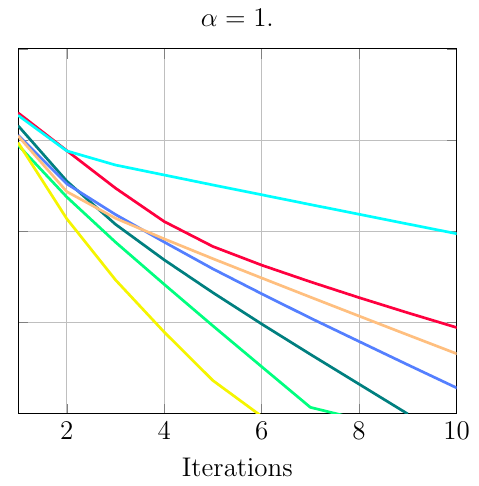}	
\raisebox{3.5ex}{
\includegraphics[width=0.15\textwidth]{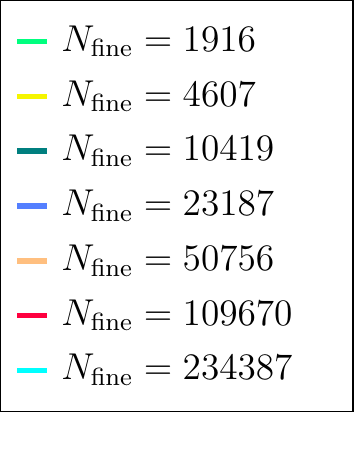}	}
\caption{$\log_{10}$ of the Euclidean norm of the residual for the multigrid method applied to the dual formulation for the face problem in Figure \ref{face_problem}. Parameters: $\mu=1$, $\lambda=\infty$, $\text{number of smoothing steps}=5$. The coarsest level has dimension $N_{\text{coarse}}=769$. Then bisection on each element is used to refine the mesh. The residuals have been computed after V-cycle. }
\label{face_problem_alpha}
\end{figure}

\begin{figure}[htbp!]
\hspace{+3cm}
\includegraphics[width=0.27\textwidth]{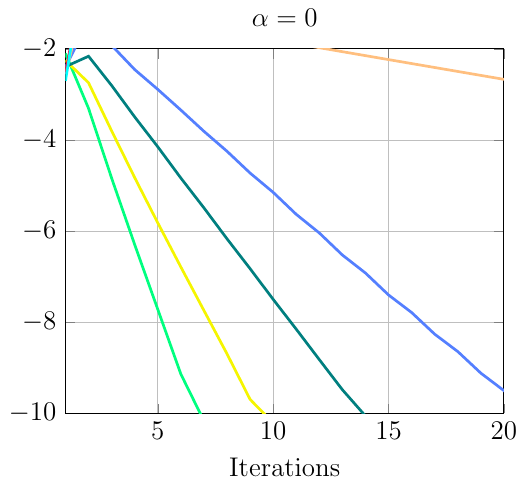}
\includegraphics[width=0.25\textwidth]{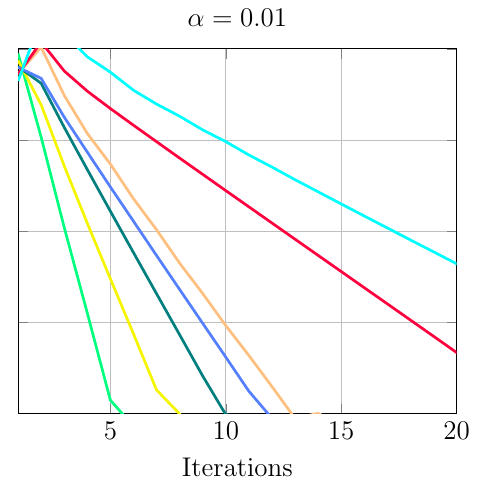}		

\hspace{+3cm}
\includegraphics[width=0.27\textwidth]{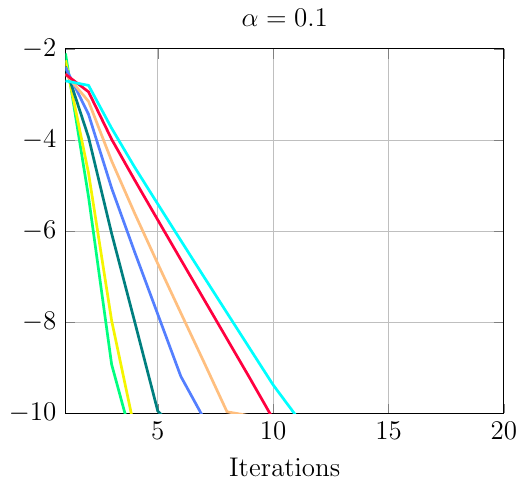}
\includegraphics[width=0.25\textwidth]{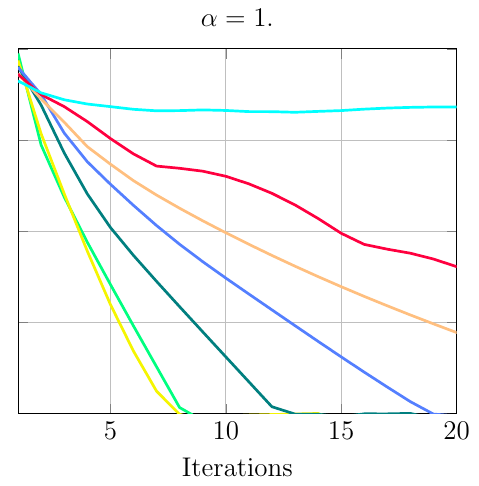}	
\raisebox{3.5ex}{
\includegraphics[width=0.15\textwidth]{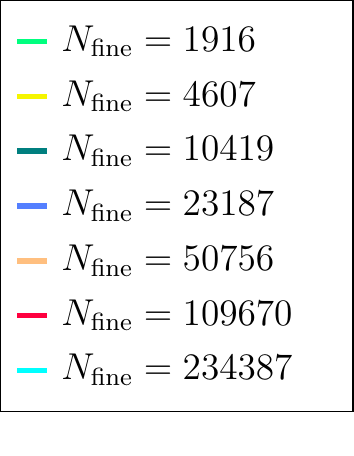}	}
\caption{$\log_{10}$ of the Euclidean norm of the residual for the 2 grids method applied to the dual formulation for the face problem in Figure \ref{face_problem}. Parameters: $\mu=1$, $\lambda=\infty$, $\text{number of smoothing steps}=5$. The coarsest level has dimension $N_{\text{coarse}}=769$. Then bisection on each element is used to refine the mesh. The residuals have been computed after V-cycle. }
\label{2grid_face_problem_alpha}
\end{figure}

\begin{figure}[htbp!]

\hspace{+2cm}
\begin{subfigure}{.35\textwidth}
\includegraphics[scale=0.2]{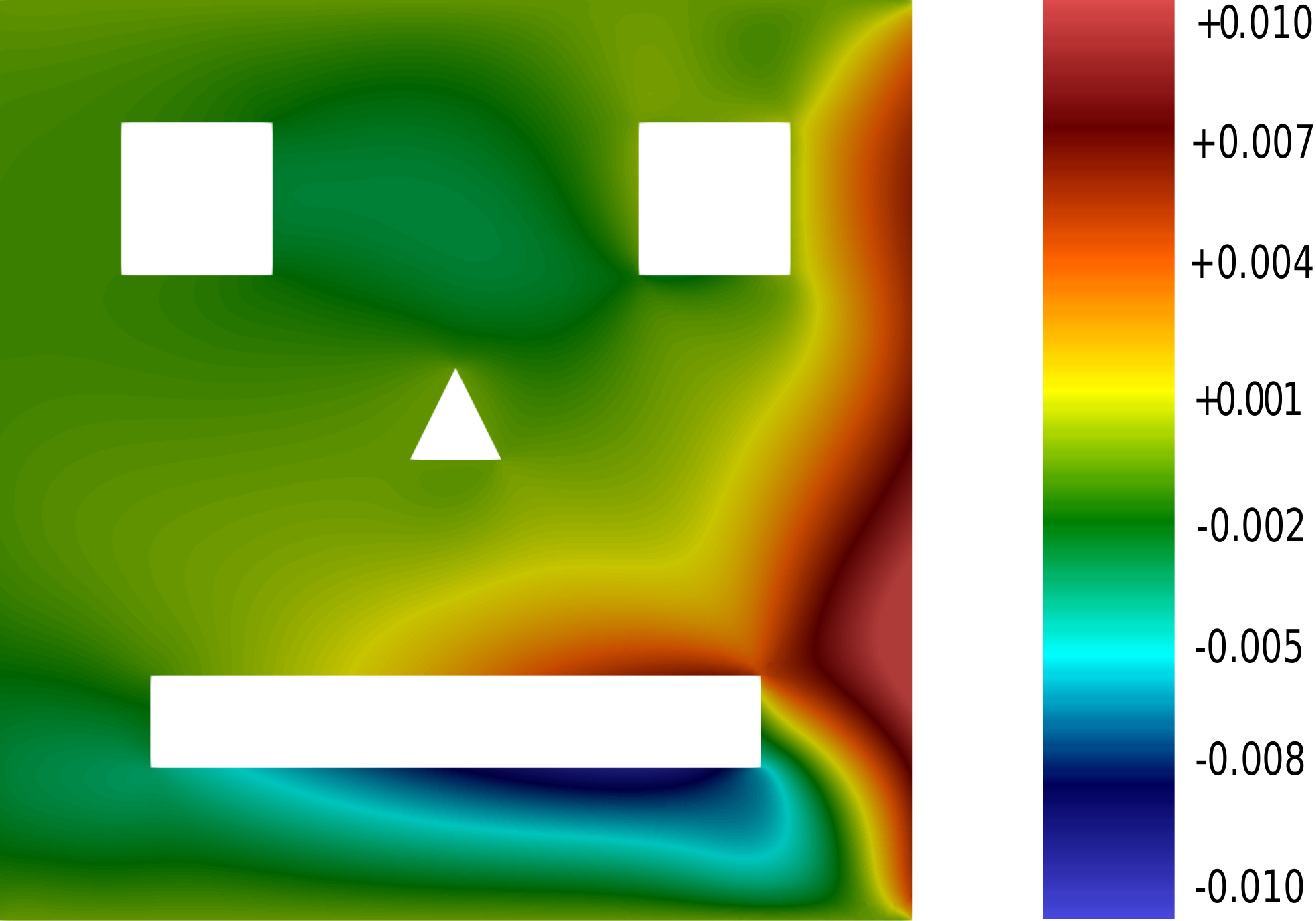}
\caption{$u_x  \qquad \qquad\qquad\qquad\qquad$}
\end{subfigure}
$\qquad$ $\qquad$
\begin{subfigure}{.31\textwidth}
\includegraphics[scale=0.2]{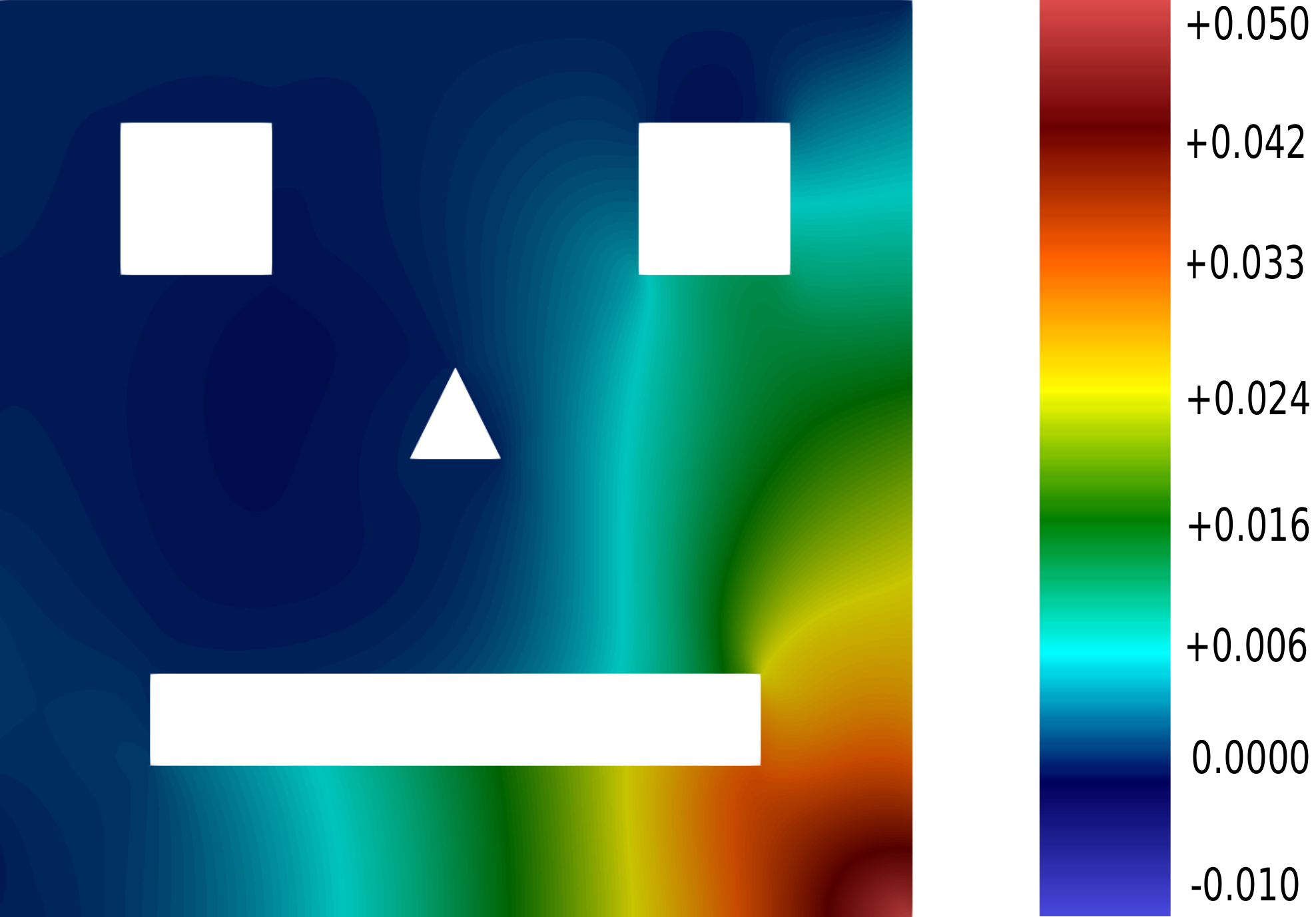}
\caption{$u_y  \qquad \qquad\qquad\qquad\quad$}
\end{subfigure}

\hspace{+2cm}
\begin{subfigure}{.35\textwidth}
\includegraphics[scale=0.2]{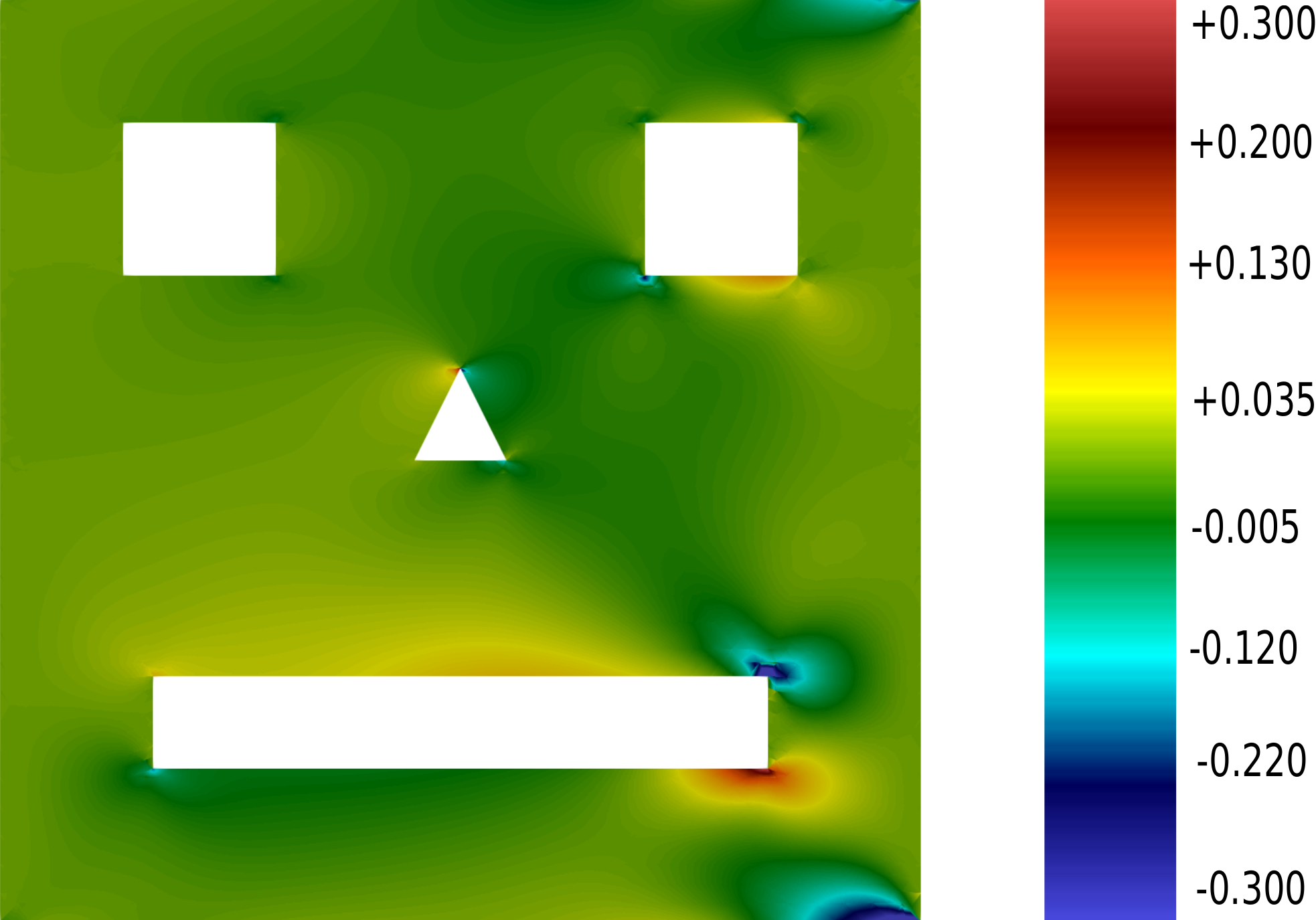}
\caption{$\sigma_{xx} \qquad \qquad\qquad\qquad\qquad$}
\end{subfigure}
$\qquad$ $\qquad$
\begin{subfigure}{.35\textwidth}
\includegraphics[scale=0.2]{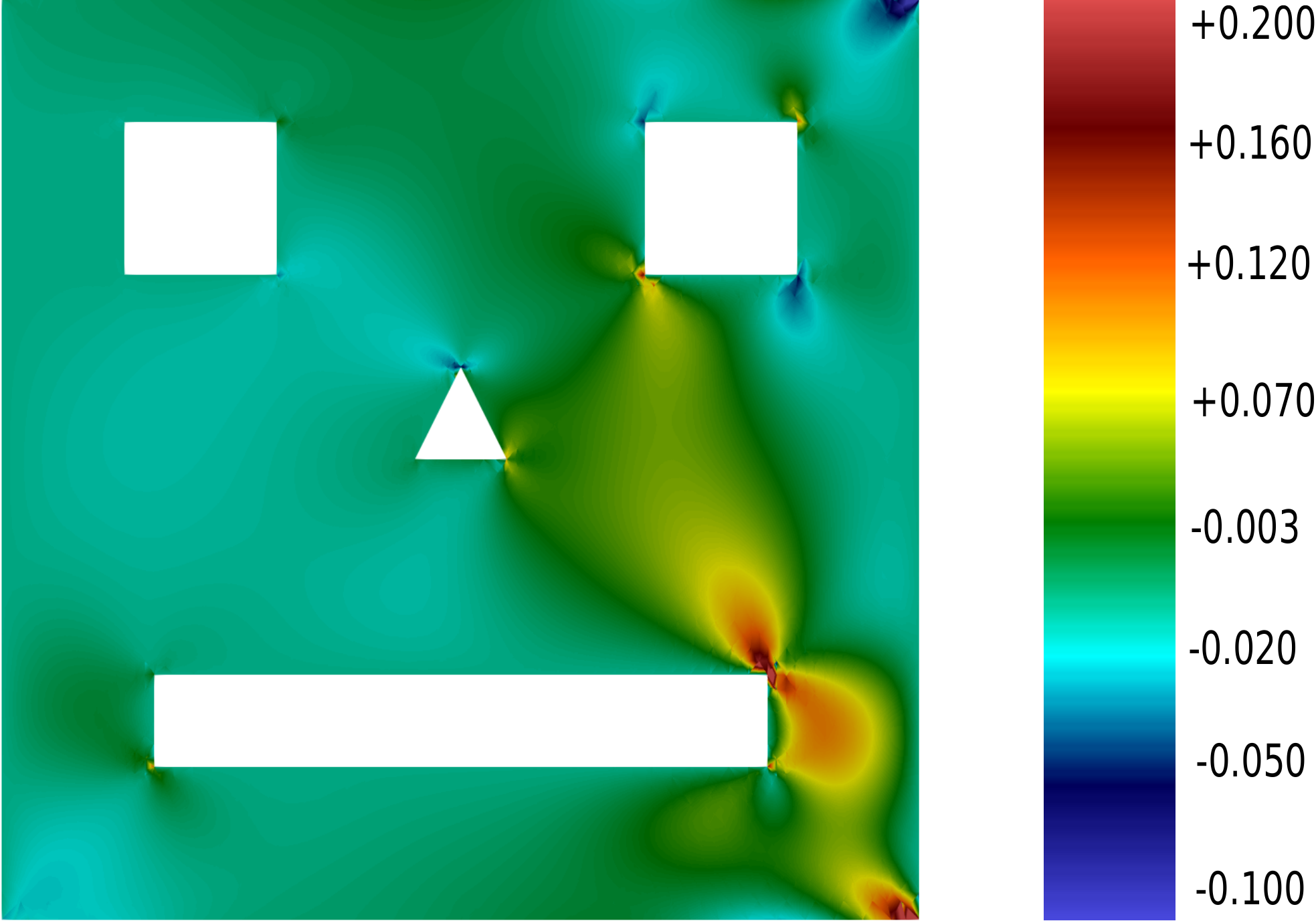}
\caption{$\sigma_{xy} \qquad \qquad\qquad\qquad\qquad$}
\end{subfigure}

\hspace{+2cm}
\begin{subfigure}{.35\textwidth}
\includegraphics[scale=0.2]{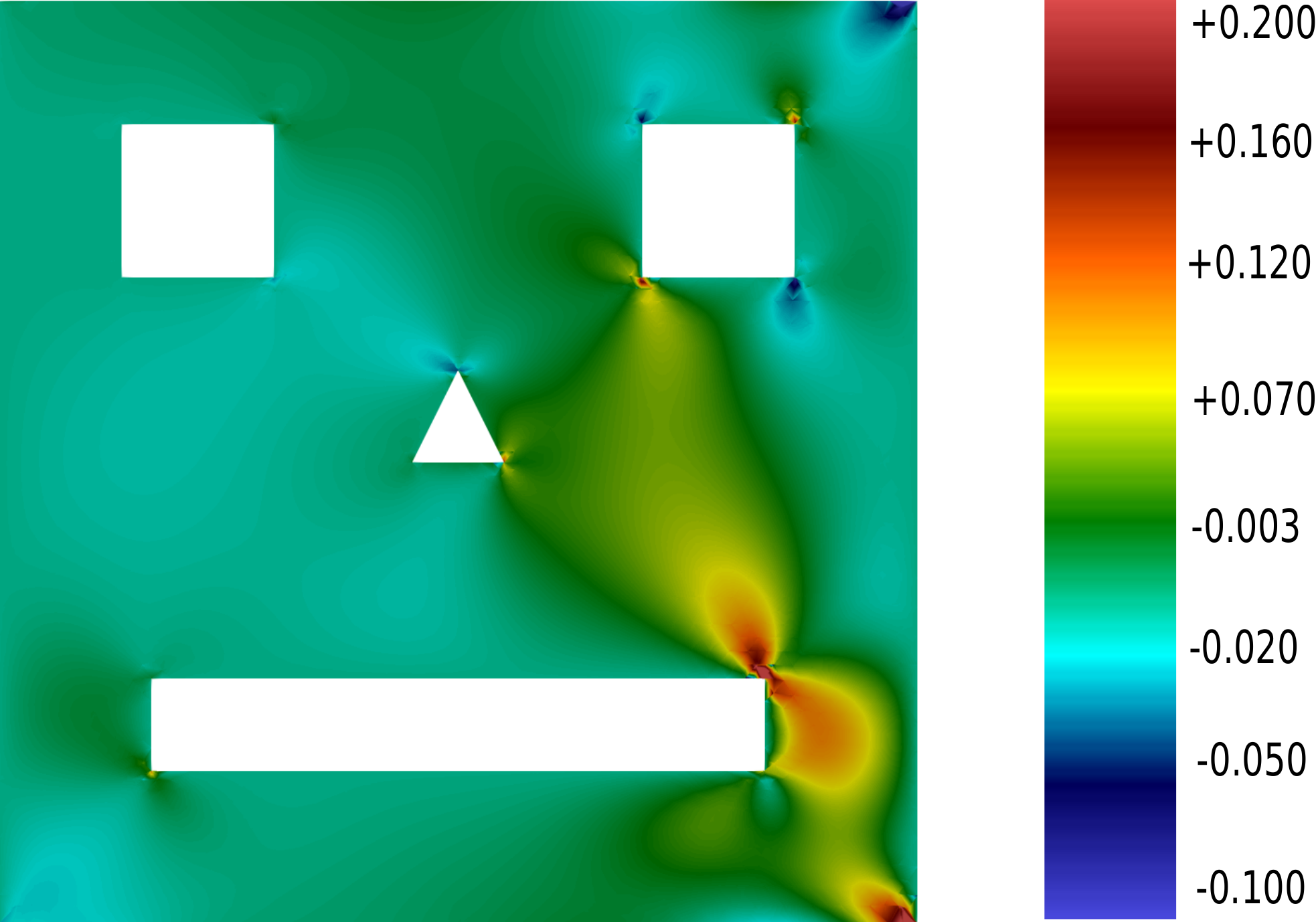}
\caption{$\sigma_{yx} \qquad \qquad\qquad\qquad\qquad$}
\end{subfigure}
$\qquad$ $\qquad$
\begin{subfigure}{.31\textwidth}
\includegraphics[scale=0.2]{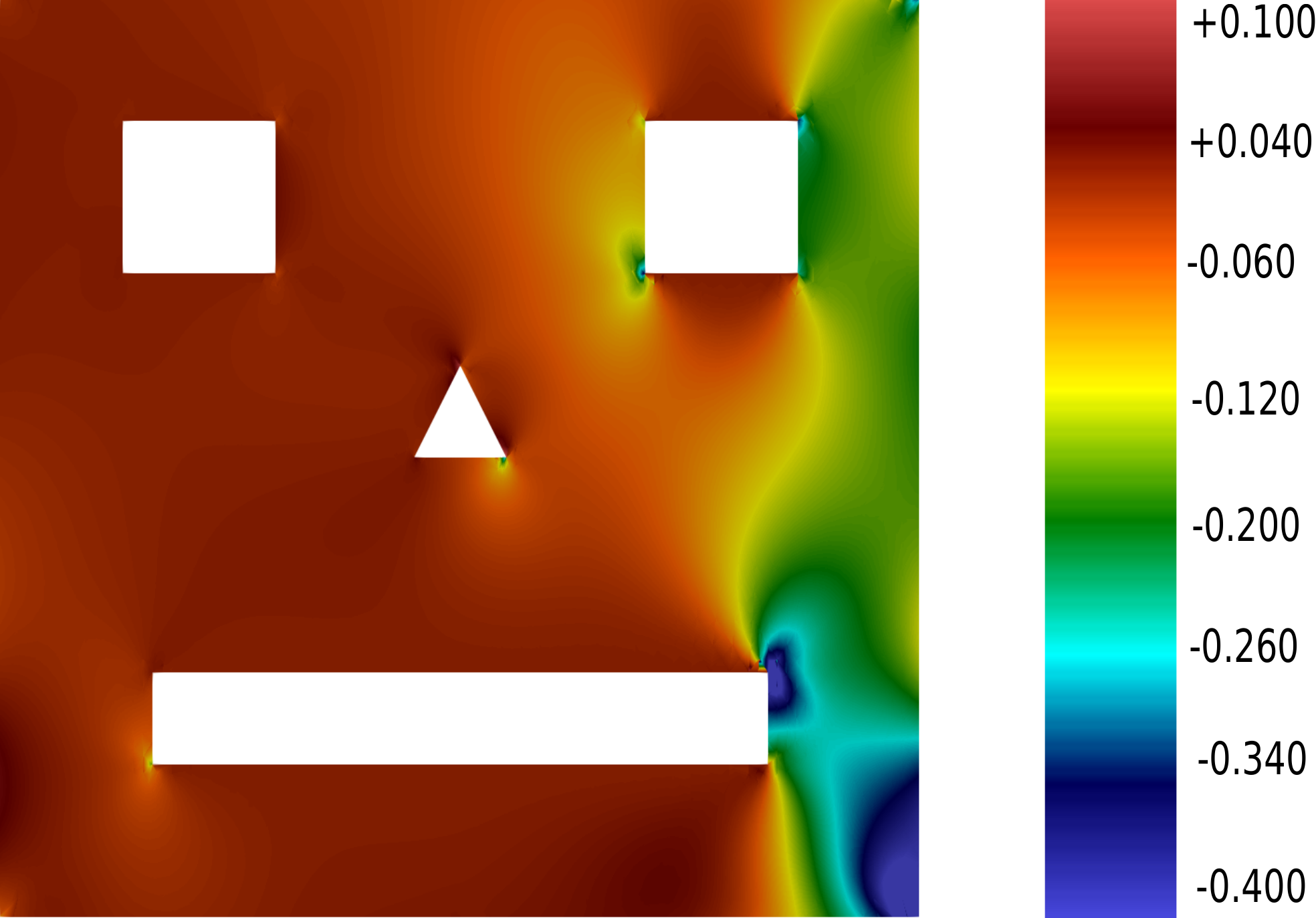}
\caption{$\sigma_{yy} \qquad \qquad\qquad\qquad\quad$}
\end{subfigure}

\hspace{+2cm}
\begin{subfigure}{.35\textwidth}
\includegraphics[scale=0.2]{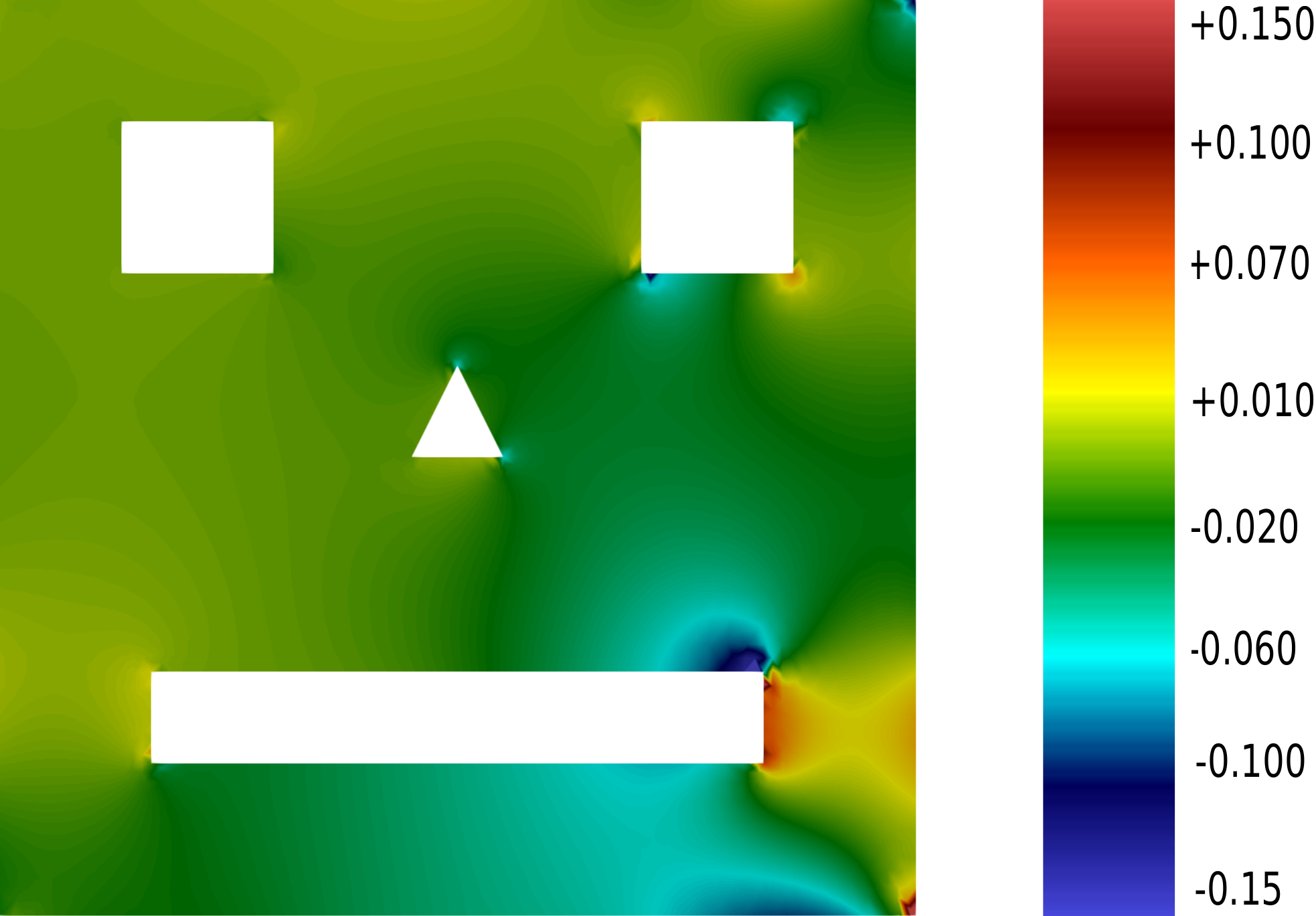}
\caption{$\rho \qquad \qquad\qquad\qquad\qquad$}
\end{subfigure}
$\qquad$ $\qquad$
\begin{subfigure}{.45\textwidth}
\includegraphics[scale=0.065]{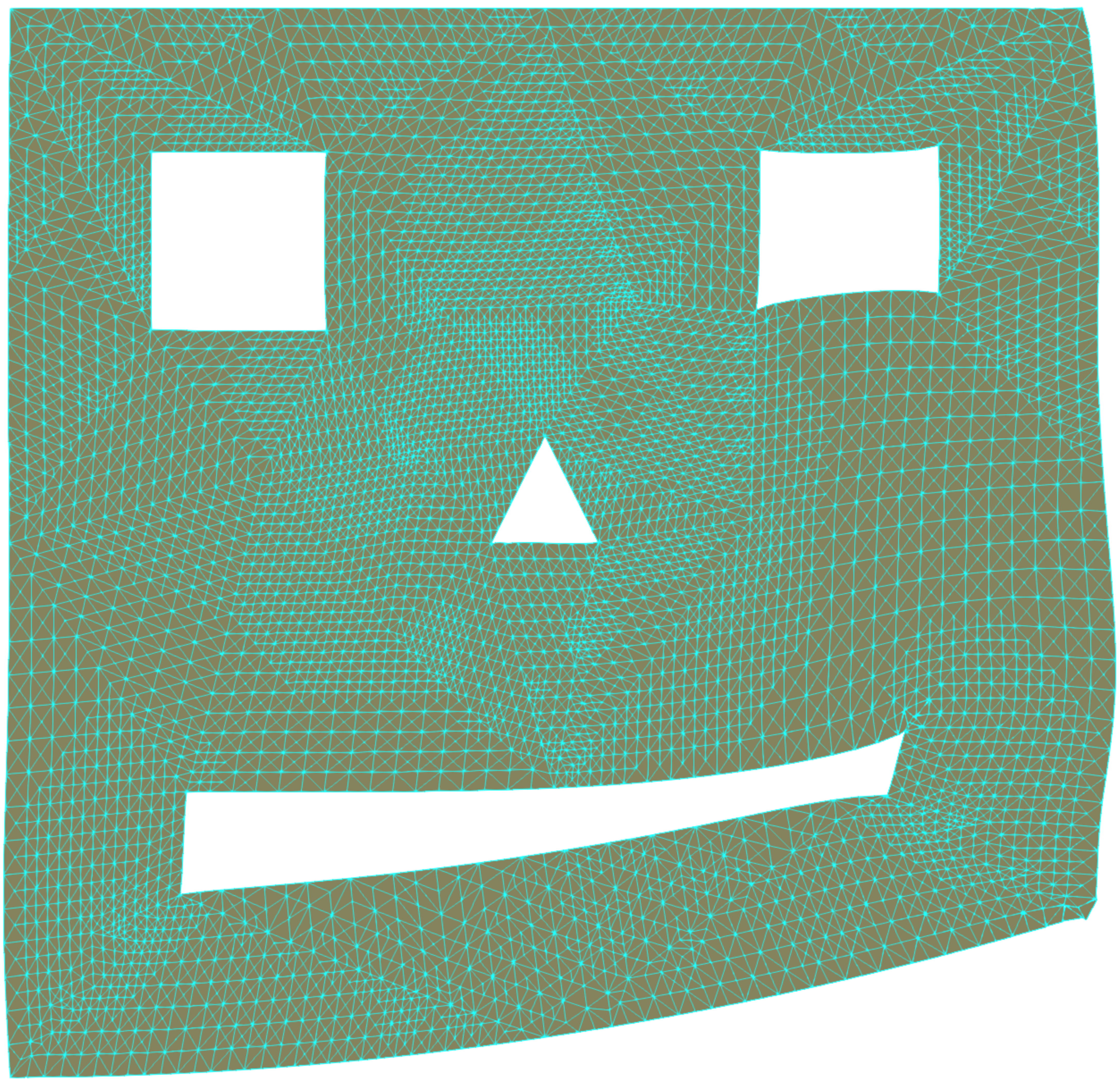}
\caption{3 times larger deformed mesh $\qquad$}
\end{subfigure}
\caption{Results for the face problem.
Parameters: $\mu=1$, $\lambda=\infty$.}
\label{face_contact_solutions}
	\end{figure}

In this section we want to examine how non-convexity of the geometry can eventually affect the multigrid convergence. Let us consider a square-shaped domain $\Omega$ with four different holes: two squares, one triangle and one rectangle, as depicted in Figure \ref{face_problem}. On the bottom edge, we enforce a quadratic displacement in the $y-$component, i.e. ${\bg_D|_{\text{bottom}}=[0,0.05 x^2]^T}$. On the triangle, we enforce zero displacement, i.e. ${\bg_D|_{\text{triangle}}=[0,0]^T}$. Everywhere else, we impose homogeneous Neumann conditions, ${\bg_N =[0,0]^T}$. The material has the following parameter: ${\mu=1}$ and ${\lambda=\infty}$. Solutions are represented in Figure \ref{face_contact_solutions}.\\
In Figure \ref{face_smoothing_alpha}, the smoother as a stand-alone solver is examined for ${\alpha=0, \:0.01,\: 0.1, \:1}$. This choice is due to the fact that, as it happens for the Cook's problem, for increasing $\alpha > 1$ the convergence behavior starts to deteriorate. Thus, we have preferred to represent more values of $\alpha$ in  the range $[0,1]$, opting for adding ${\alpha=0.01,\:0.1}$. The plots in Figure \ref{face_problem_alpha} are obtained using a multigrid method with a coarse mesh of ${N_{\text{coarse}}=769}$ dofs which we refine, with a bisection algorithm, up to ${N_{\text{fine}}=234387}$ dofs. For each level, we do 5 pre-smoothing steps and 5 post-smoothing steps. On the coarsest level, we solve exactly. We have repeated the experiments for the same values of  ${\alpha=0, \:0.01,\: 0.1, \:1}$. Results are represented in Figure \ref{face_problem_alpha}. In contrast to Cook's membrane problem, the parameter $\alpha$ has to be chosen carefully. If $\alpha \leq 0.01$ results are always good, but for larger values of $\alpha$, the rate of convergence is no more optimal and depends on the size of the problem. For too large $\alpha$, the method does not even converge. By inspecting the Figure \ref{2grid_face_problem_alpha} for aggressive coarsening, it is clear that again the full Dirichlet case $\alpha=0$ does not give rise to a convergent method. On the other hand, the optimal behavior is reached for $\alpha=0.1$. Indeed the rate of convergence is independent of the dimension of the problem. This means that, given a coarse level, independently of the fine level considered, the Robin boundary conditions with a proper value of $\alpha$ can damp all the frequency components of the error in between.
By moving far away from this value, convergence deteriorates or is lost.


\section{Conclusion}
In this paper, a smoother for the dual linear elasticity problem is introduced. Since the stress belongs to $\bH_{\text{div}}$, we take advantage of the patch-smoother by Arnold-Falk-Winther, that we also extend to degrees of freedom related to the Lagrange multipliers, i.e. the displacement and the rotation. Indeed, in the incompressible limit, the stress block of the corresponding saddle point system is only semi-positive definite and thus not invertible. However different possibilities on the boundary conditions for the local stress can be examined. The full Neumann approach gives rise to a divergent smoother, while the full Dirichlet approach makes the method convergent. However, in the latter case, the violation of the constraints outside the patch suggests the damping of the local stress correction on the boundary patch. The damped system can be interpreted as a system where local Robin boundary conditions of parameter $\alpha$ are enforced. For $\alpha=0$, the Dirichlet system is recovered. In the numerical experiments, it is shown that for the multigrid case, $\alpha$ can be chosen to be zero with no problem and too large values could affect badly the convergence. However, for the aggressive coarsening case, we have shown that a proper tuning of $\alpha$ is fundamental and that $\alpha=0$ can make the method diverge. In particular, for future generalizations to non-linear cases, like the Signorini problem, we can also expect the necessity of choosing carefully the parameter $\alpha$.
\section*{Acknowledgments}
The authors would like to thank the Swiss National Science Foundation for their support through the project and the Deutsche Forschungsgemeinschaft (DFG) for their support in the SPP 1962 ``
Stress-Based Methods for Variational Inequalities in Solid Mechanics: Finite Element Discretization and Solution by Hierarchical Optimization [186407]''.

\end{sloppypar}

\bibliographystyle{abbrv}
\bibliography{references.bib}
\end{document}